\documentclass{article}

\usepackage[textsize=footnotesize]{todonotes}
\usepackage{amsmath}
\usepackage{amsthm}
\usepackage{amsfonts}
\usepackage{algorithm}
\usepackage{algorithmic}
\usepackage{array}
\usepackage[affil-it]{authblk}
\usepackage{bm}
\usepackage{bbm}
\usepackage{booktabs}
\usepackage{color}
\usepackage[OT1]{fontenc}
\usepackage{graphicx}
\usepackage{float}
\usepackage{bigints}

\usepackage{multirow}
\usepackage{subcaption}

\usepackage{mathrsfs}
\usepackage[bbgreekl]{mathbbol}
\usepackage{mathabx}
\usepackage{mathtools}
\usepackage{natbib}
\usepackage[section]{placeins}
\usepackage{verbatim}
\usepackage{enumitem}

\usepackage{mdframed}

\usepackage{xr}
\usepackage{xcolor}

\usepackage{hyperref}
\usepackage{float}

\DeclareSymbolFontAlphabet{\mathbbm}{bbold}
\DeclareSymbolFontAlphabet{\mathbb}{AMSb}%

\usepackage{setspace}
\newtheorem{theorem}{Theorem}[section]
\newtheorem{lemma}{Lemma}[section]
\newtheorem{proposition}{Proposition}[section]

\newtheorem{remark}{Remark}[section]

\numberwithin{equation}{section}
\numberwithin{theorem}{section}
\numberwithin{lemma}{section}
\numberwithin{proposition}{section}
\numberwithin{corollary}{section}
\numberwithin{definition}{section}
\numberwithin{cons}{section}
\numberwithin{remark}{section}
\numberwithin{exa}{section}
\numberwithin{table}{section}
\numberwithin{figure}{section}

\newcommand{\nonsig}{D}

\newcommand{\bX}{\mathbf{X}}
\newcommand{\bA}{\mathbf{A}}
\newcommand{\bB}{\mathbf{B}}

\newcommand{\bW}{\mathbf{W}}

\newcommand{\calA}{\mathcal{A}}
\newcommand{\calB}{\mathcal{B}}

\newcommand{\calG}{\mathcal{G}}
\newcommand{\calH}{\mathcal{H}}
\newcommand{\calI}{\mathcal{I}}

\newcommand{\frakM}{\mathfrak{M}}
\newcommand{\frakN}{\mathfrak{N}}

\newcommand{\supp}{\operatorname{supp}}

\begin{document}
\title{Analysis of The Limiting Spectral Distribution of Large Random Matrices of The Mar\v{c}enko-Pastur Type}
\author{Haoran Li\footnote{Department of Mathematics and Statistics, Auburn University, 221 Parker Hall, Auburn, AL, 36849, hzl0152@auburn.edu}}
\affil{Auburn University}
\maketitle

\allowdisplaybreaks

\begin{abstract}
\setlength{\baselineskip}{1.665em}
Consider the random matrix \(\bW_n = \bB_n + n^{-1}\bX_n^*\bA_n\bX_n\), where \(\bA_n\) and \(\bB_n\) are Hermitian matrices of dimensions \(p \times p\) and \(n \times n\), respectively, and \(\bX_n\) is a \(p \times n\) random matrix with independent and identically distributed entries of mean 0 and variance 1. Assume that \(p\) and \(n\) grow to infinity proportionally, and that the spectral measures of \(\bA_n\) and \(\bB_n\) converge as \(p, n \to \infty\) towards two probability measures \(\calA\) and \(\calB\). Building on the groundbreaking work of \cite{marchenko1967distribution}, which demonstrated that the empirical spectral distribution of \(\bW_n\) converges towards a probability measure \(F\) characterized by its Stieltjes transform, this paper investigates the properties of \(F\) when \(\calB\) is a general measure. We show that \(F\) has an analytic density at the region near where the Stieltjes transform of $\calB$ is bounded. The density closely resembles \(C\sqrt{|x - x_0|}\) near certain edge points \(x_0\) of its support for a wide class of \(\calA\) and \(\calB\). We provide a complete characterization of the support of \(F\). Moreover, we show that \(F\) can exhibit discontinuities at points where \(\calB\) is discontinuous.
\\
\noindent {\bf Keywords}: Random Matrix Theory,  Mar\v{c}enko-Pastur distribution, Stieltjes transform.
\end{abstract}


\section{Introduction}\label{sec:introduction}

For $n=1,2,\dots$ and $p=p(n)$, let $\bX_n = (X_{ij}^n)$ is a $p\times n$ matrix of complex random variables that are independent for each $n$,  identically distributed for all $n$, and satisfy $\mathbb{E} X_{11} = 0$ and  $\mathbb{E}\left|X_{11}\right|^2=1$. Let \(\bA_n\) be a diagonal matrix and \(\bB_n\) be a Hermitian random matrix, independent of \(\bX_n\). In Random Matrix Theory, a class of matrices of the form
\[
\bW_n = \bB_n + n^{-1} \bX^*_n \bA_n \bX_n
\]
is of significant interest, where \(\bX^*_n\) stands for the conjugate transpose of \(\bX_n\). The behavior of the eigenvalues of this matrix when \(n\) and \(p\) are both large and comparable was first studied in the groundbreaking work of \cite{marchenko1967distribution}. Since then, many researchers made significant contributions to extending the results and relaxing the assumptions. Notable works include \cite{grenander1977spectral, wachter1978strong, jonsson1982some, yin1983limit, yin1986limiting, silverstein1995empirical}.
The behavior of the eigenvalues is expressed in terms of a limit theorem, as \(n \to \infty\) while \(p = p(n)\) with \(p/n \to \gamma > 0\), on the empirical distribution function (e.d.f.) \(F^{\bW_n}\) of the eigenvalues of \(\bW_n\). Specifically, \(F^{\bW_n}(x)\) represents the proportion of eigenvalues of \(\bW_n\) that are less than or equal to \(x\). It has been concluded that \(F^{\bW_n}\) converges in some sense to a non-random distribution \(F\). The type of convergence varies under different assumptions on \(\bA_n\) and \(\bB_n\). The aim of the present paper is to derive analytical properties of \(F\) when $\bA_n$ and $\bB_n$ are general.

To introduce the limiting behavior of \(F^{\bW_n}\), we start with the Stieltjes transform, which plays a critical role in Random Matrix Theory. Recall that the Stieltjes transform of a measure \(\mu\) on the real line is defined as the function
\[
m_\mu(z) = \int \frac{1}{x - z} \mu(dx), \quad z \in \mathbb{C}^+,
\]
where \(\mathbb{C}^+ = \{x + iy: x \in \mathbb{R},~ y > 0\}\). Note that \(m_\mu\) is analytic on \(\mathbb{C}^+\) and maps into \(\mathbb{C}^+\). We can retrieve \(\mu\) using the inversion formula
\[
\mu\{[a, b]\} = \frac{1}{\pi} \lim_{\eta \rightarrow 0^{+}} \int_a^b \Im m_\mu(\xi + i \eta) \, d\xi,
\]
where \(a\) and \(b\) are continuity points of \(\mu\). Following \cite{silverstein1995empirical}, assume:
\begin{itemize}
    \item[\textbf{C1}] \(\bA_n = \text{diag}(u_1^n, \dots, u_p^n)\), where \(u_i^n \in \mathbb{R}\) is a deterministic matrix and the e.d.f. of \((u_1^n, \dots, u_p^n)\) converges weakly to a probability distribution function \(\calA\) as \(n \to \infty\). Without loss of generality, assume that \(\calA\) is not the Dirac measure at 0.
    \item[\textbf{C2}] The matrix $\bB_n$ is Hermitian and the e.d.f. of the eigenvalues of \(\bB_n\) converges weakly to \(\calB\) almost surely, where \(\calB\) is a non-random distribution function.
\end{itemize}
Then, almost surely, as \(n, p \to \infty\) and \(p/n \to \gamma \in (0,\infty)\), \(F^{\bW_n}\) converges weakly to a non-random distribution \(F\) whose Stieltjes transform \(m(z)\), for \(z \in \mathbb{C}^+\), is determined as the unique solution in \(\mathbb{C}^+\) to the equation
\begin{equation}\label{eq:main_eq}
m(z) = m_{\calB}\left(z - \gamma \int \frac{u \, d\calA(u)}{1 + m(z) u}\right) = \int \frac{d\calB(\tau)}{\tau - z + \gamma \displaystyle\int \frac{u \, d\calA(u)}{1 + m(z) u}}.
\end{equation}
The convergence and the existence and uniqueness of \(m(z)\) for any \(z \in \mathbb{C}^+\) are shown in \cite{marchenko1967distribution} and \cite{silverstein1995empirical}.

Unfortunately, no explicit expression for $F$ exists except when $\calA$ and $\calB$ are simple. However, quite a bit of properties of $F$ can be inferred from the master equation \eqref{eq:main_eq} of \(m(z)\). \cite{silverstein1992signal} and \cite{silverstein1995analysis} pioneered research in this direction. They considered the sample covariance matrix, which corresponds to the simplified case of \(\bW_n\) when \(\bB_n = 0\). They proved that the density function of $F$ exists away from zero and further characterized its support. The density is real analytic wherever it is positive and exhibits only square-root singularities at the boundaries (edges) of its support. This square-root behavior is a fundamental observation in Random Matrix Theory, as it is closely linked to the Tracy-Widom fluctuations of eigenvalues near these edges. It has been used as a critical condition in \cite{elkaroui2007tracy, erdos2013local, bloemendal2014isotropic, lee2016tracy, knowles2017anisotropic}. Additionally, the universality and rigidity of eigenvalues of sample covariance matrices are intimately connected to these properties, as explored in \citet{bai1998no, pillai2014universality, bao2015universality}. On the application side, these results have been leveraged in \citet{mestre2008improved, mestre2008asymptotic} to develop consistent estimators for linear functionals of the eigenvalues of $\bA_n$.

\citet{couillet2014analysis} extended results on sample covariance matrices to random matrices of the separable covariance type. A main contribution of their work is the design of more fundamental analytical tools for dealing with the master equation of the Stieltjes transform. The technical arguments in this article closely follow their approach. \citet{dozier2007analysis} investigated the limiting spectral distribution of information-plus-noise type matrices. \citet{ajanki2017singularities} studied solutions to a class of quadratic vector equations that naturally arise in the analysis of Wigner-type matrices with general variance structures. They established the presence of square-root singularities at the edges of the support and cubic-root cusps in the interior. Similar results for random Gram matrices were obtained by \citet{alt2017local} and \cite{Alt2017singularities}. \citet{Alt2020} analyzed a class of correlated random matrices of the form $B+ H$ where $B$ is Hermitian and deterministic, and $H$ is Hermitian with correlated entries. They demonstrated the square-root behavior at the support edges and the Tracy-Widom fluctuations of eigenvalues near these edges. Other notable works include \citet{ajanki2019stability}, \citet{ajanki2017universality}, \citet{alt2020dyson}, and \citet{alt2019location}.

The current work analyzes the limiting spectral distribution of the sum of a sample covariance matrix and a general independent random matrix with a convergent spectral measure. It extends the results of \citet{silverstein1995analysis} to the case where $\bB_n$ is nonzero. Notably, except when $\bB_n$ is a multiple of the identity matrix (i.e., when $\calB$ is a Dirac measure), the problem becomes nontrivial. The primary application of this framework is the spectral analysis of regularized sample covariance matrices. While the sample covariance matrix is inconsistent in high-dimensional settings, regularization by adding a matrix $\bB_n$ has gained attention for mitigating bias in high-dimensional inference. For example, inspired by ridge regression, \citet{li2020adaptable} and \citet{li2020high} studied ridge-regularized sample covariance matrices where $\bB_n$ is a multiple of the identity matrix. \citet{elkaroui2011geometric} explored more flexible regularization, allowing $\bB_n$ to take general forms, and established the concentration of a quadratic form involving the regularized sample covariance matrix. These studies have laid the groundwork for incorporating flexible regularization into statistical inference. Understanding the limiting spectral distribution is therefore crucial for analyzing the asymptotic properties of regularized sample covariance matrices.

The rest of the paper is organized as follows. In Section \ref{sec:preliminary}, we present some preliminary results. In Section \ref{sec:existence}, we show the existence of a density function of \(F\) in the region near where the Stieltjes transform of $\calB$ is bounded. The density function is shown to be analytic when positive.   We further discuss the mass of \(F\) at discontinuity points of \(\calB\). In Section \ref{sec:support}, we characterize the support of \(F\) and present a practical procedure to determine the support. In Section \ref{sec:boundary_behavior}, we show that the behavior of the density function near an edge point of its support.

\section{Preliminary results}\label{sec:preliminary}
Throughout the paper, we assume that \(\calB\) is not a Dirac measure. If \(\calB\) is the Dirac measure at a point \(a \in \mathbb{R}\), after a shift by \(-a\), \(F\) reduces to the limiting distribution of \(n^{-1}\bX^* \bA_n \bX_n\), whose properties are studied in \cite{silverstein1995analysis}. In addition to \textbf{C1} and \textbf{C2}, we impose the following condition on \(\calA\):
\begin{itemize}
    \item[\textbf{C3}] The second moment of \(\calA\) exists.
\end{itemize}

For any \(x \in \mathbb{R}\), we define its \(\varepsilon\)-neighborhood in \(\mathbb{C}^+\) as 
$\frakN(x, \varepsilon) = \{z : |z - x| \leq \varepsilon \text{ and } z \in \mathbb{C}^+\}$. 
Define 
\begin{equation}\label{eq:def_D}
D = \Big\{ x \in \mathbb{R} : \text{there exists } \varepsilon > 0 \text{ such that } m_{\calB}(z) \text{ is bounded in } z \in \frakN(x, \varepsilon) \Big\}.
\end{equation}
Clearly, \(D\) is open.   We have the following result. 
\begin{lemma}\label{lemma:element_of_D}
If $\calB(x)$ is analytic at $x$, we have $x \in \nonsig$.
\end{lemma}
\noindent The lemma directly follows from the residue theorem, and the proof is omitted. 
It is clear that if \(x\) is outside the support of \(\calB\), then \(x \in D\). The interior points of the support of $\calB$ are also in $D$ if the density is smooth.
On the other hand, if the density of \(\calB\) has a pole at \(x\), then \(x \in D^\complement\). Discontinuity points of \(\calB\) also lie in \(D^\complement\). We provide three representative examples:
\begin{itemize}
\item[(i)] If $\calB$ is a discrete measure with masses on $\{\tau_1,\tau_2,\dots, \tau_k\}$, it is straightforward that $D = \mathbb{R}\setminus\{\tau_1,\tau_2,\dots, \tau_k\}$.
\item[(ii)] If $\calB$ is the semicircular law with the radius parameter $R^2$, whose pdf is given as  $f_{\calB}(x) = ({2}/{\pi R^2})\sqrt{R^2 -x^2}$, $-R \leq x\leq R$, the Stieltjes transform is given as $m_{\calB}(z) = ({-2}/{R^2})(z -\sqrt{z^2 - R^2})$, $z\in\mathbb{C}^+$ (see, for example, Lemma 2.11 of \cite{bai2010spectral}). Here, the complex square root is taken to have positive imaginary part.
It follows that $D = \mathbb{R}$. It is worth pointing out that $f_{\calB}(x)$ is not analytic at $\pm R$. Lemma \ref{lemma:element_of_D} is not a necessary condition. 
\item[(iii)] If $\calB$ is the Mar\v{c}enko-Pastur distribution with parameter $\lambda>1$, $\calB$ is a mixture of a continuous measure and a point-mass at $0$. Following \cite{silverstein1995analysis}, $\lim_{z\to x}m_{\calB}(z)$ exists for any $x\in\mathbb{R}\neq 0$, implying $D=\mathbb{R}\setminus\{0\}$. Indeed, the density of $\calB$ is analytic where it is positive. 
\end{itemize}

We define some auxiliary objects that play a central role in the sequel. For $z\in\mathbb{C}^+$, call its conjugate as $z^*$. For any $z, z_1,z_2\in\mathbb{C}^+$, define 
\begin{align*}
&\displaystyle\alpha(z_1, z_2) =\bigintssss \frac{d\calB(\tau)}{\left(\tau - z_1 + \displaystyle\gamma\int\frac{ud\calA(u)}{1+um(z_1)}\right)\left(\tau - z_2 + \displaystyle\gamma\int \frac{ud\calA(u)}{1+um(z_2)}\right) }, \\[7pt]
&\beta(z_1, z_2) = \gamma \int \frac{u^2d\calA(u)}{(1+um(z_1))(1+um(z_2))},\\[7pt]
&\alpha(z, z^*) = \bigintssss \frac{d\calB(\tau)}{\left|\tau - z + \displaystyle\gamma\int\frac{ud\calA(u)}{1+um(z)}\right|^2}, \\[7pt]
&\beta(z, z^*) = \gamma\int\frac{u^2 d\calA(u)}{|1+um(z)|^2}.
\end{align*}
Note that  
\[\left|\tau - z + \gamma\int\frac{ud\calA(u)}{1+um(z)}\right|\geq \Im z - \Im \gamma\int\frac{ud\calA(u)}{1+um(z)} \geq \Im z.\] 
These integrands above are integrable. 
The imaginary part of both sides of \eqref{eq:main_eq} can be concisely written as 
\begin{align*}
 \Im m(z) &= \left(\Im z+ \gamma\Im m(z)\int\frac{u^2d\calA(u)}{|1+um(z)|^2 } \right)\int \frac{d\calB(\tau)}{\Big|\tau -z + \displaystyle\gamma\int\frac{ud\calA(u)}{1+um(z)} \Big|^2}\\
&= \Big(\Im z + \beta(z,z^*) \Im m(z)\Big) \alpha(z,z^*).
\end{align*}
It follows that for any $z\in\mathbb{C}^+$,
\begin{equation}\label{eq:ineq_core}
\Big(1 -\alpha(z,z^*)\beta(z,z^*)\Big)=\frac{\Im z }{\Im m(z)} \alpha(z,z^*)>0.
\end{equation}

\section{Density function and mass}\label{sec:existence}
In this section, we show the following results.  
\begin{theorem}\label{thm:main}
For any $x\in \nonsig$, we have  $\lim_{z\in\mathbb{C}^+ \to x} m(z)\equiv \underline{m}(x)$  exists. The function $\underline{m}(x)$ is continuous on $\nonsig$, and $F$ has a continuous derivative $f$ on $D$ given by 
\begin{equation}\label{eq:expression_density}
f(x) = \frac{1}{\pi} \Im \underline{m}(x).
\end{equation}
Moreover, if $x\in\nonsig$ is such that $f(x)>0$, $f(x)$ is analytic at $x$. 
\end{theorem}
\begin{theorem}\label{thm:existence_of_masses}
Denote the mass of $\calA$ at $0$ by $\calA(\{0\})$. 
For any $b\in\mathbb{R}$, the necessary and sufficient condition for $F$ to have  positive mass at $b$ is $\calB(\{b\}) - \gamma (1- \calA(\{0\}))>0$. In this case, the mass is $F(\{b\}) = \calB(\{b\}) - \gamma (1- \calA(\{0\}))$.
\end{theorem}
Theorem \ref{thm:main} indicates that $F$ has a smooth density function where $\calB$ is smooth. The result is a generalization of Theorem 1.1 of \cite{silverstein1995analysis}. In their work, $\bB_n =0$ (or $\calB$ is the Dirac measure at $0$) is considered and therefore $D = \mathbb{R}\setminus\{0\}$. Similar results are obtained for random matrices of the separable covariance type in \cite{couillet2014analysis} and the information-plus-noise type matrices in \cite{dozier2007analysis}. 
It is worth pointing out that the condition that $x\in \nonsig$ is not necessary for the existence of $\underline{m}(x) = \lim_{z\in\mathbb{C}^+\to x} m(z)$. For example, it can happen that a discontinuity point $b$ of $\calB$ is outside the support of $F$, then $\lim_{z\in\mathbb{C}^+\to b} m(z) = \int (\tau - b)^{-1}dF(\tau)$ exists, but $x\in\nonsig^\complement$.

Theorem \ref{thm:existence_of_masses} indicates that $F$ can only exhibit discontinuity at points where $\calB$ places large mass. To explain the result, by Proposition 2.2 of \cite{couillet2014analysis},  \(1- \gamma(1-\calA(\{0\}))\) is the mass that the asymptotic spectral distribution of $n^{-1}\bX^*_n\bA_n \bX_n$ places at $0$. We can take $\operatorname{rank}(n^{-1}\bX_n^*\bA_n\bX_n)\simeq n \gamma(1-\calA(\{0\}))$ for large $n$.  If \(\bB_n\) has an eigenvalue \(b\) of multiplicity \(n\calB(\{b\})\), after adding \(n^{-1}\bX_n^*\bA_n\bX_n\), approximately  \(\operatorname{rank}(n^{-1}\bX_n^* \bA_n \bX_n)\) of these eigenvalues will be shifted elsewhere, given that \(\calB(\{b\}) - \gamma (1- \calA(\{0\})) > 0\). If \(\calB(\{b\}) - \gamma (1- \calA(\{0\})) \leq  0\), all these eigenvalues will be shifted. 
It is worth mentioning that  if $\calB(\{b\}) - \gamma (1- \calA(\{0\})\leq 0$, the fact that $F(\{b\}) = 0$ gives no information on whether \(b\) is outside the support of \(F\). The characterization of the support of \(F\) is presented in Section \ref{sec:support}.

The rest of this section is devoted to the proof of the two theorems.  Some techniques used here are motivated by those used in the proof of Theorem 3.1 and Proposition 2.2 of \cite{couillet2014analysis}. 

\subsection{Proof of Theorem \ref{thm:main}}\label{sec:existence_subsec:proof_density}
First of all, since $m(z)$ is a Stieltjes transform of a probability measure, it is clear that $m(z)$ is holomorphic on $\mathbb{C}^+$. Due to Theorem \ref{thm:Silverstein1995_2} in Appendix \ref{appendix:lemmas}, if we can show the existence of $\underline{m}(x)$, the continuity of $\underline{m}(x)$ directly follows. Moreover, the expression of the density function Eq. \eqref{eq:expression_density} is the conclusion of Theorem \ref{thm:Silverstein1995_1}. We only need to show the existence of $\underline{m}(x)$ and the analyticity of $\Im \underline{m}(x)$ when positive. We establish a series of lemmas to achieve this.

\begin{lemma}\label{lemma:boundedness}
Consider any fixed $x\in \nonsig$. There exists $\varepsilon>0$ such that $|m(z)|$ is bounded on  $\frakN(x,\varepsilon)$. 
\end{lemma}
\begin{proof}
Suppose not. Then, we can find a sequence $z_n\in\mathbb{C}^+$, $n=1,2,\dots$ such that $z_n \to x $ and $|m(z_n)|\to\infty$ as $n\to\infty$. Using Eq. \eqref{eq:main_eq}  and Cauchy-Schwarz inequality,
\begin{equation}\label{eq:proof_lemma_boundedness_eq1}
|m(z_n)| \leq \int\left|\frac{d\calB(\tau)}{\tau - z_n + \gamma\displaystyle\int\frac{ud\calA(u)}{1+um(z_n)}   }\right| \leq \left(\int\frac{d\calB(\tau)}{\left|\tau - z_n + \gamma\displaystyle\int\frac{ud\calA(u)}{1+um(z_n)}\right|^2}\right)^{1/2} = (\alpha(z_n,z^*_n))^{1/2}.
\end{equation} 
It follows that $\alpha(z_n,z_n^*)\to\infty$, as $n\to\infty$. By Eq. \eqref{eq:ineq_core},
\[\gamma \left|\int\frac{ud\calA(u)}{1+um(z_n)}  \right|^2 \leq \beta(z_n, z_n^*) <\frac{1}{\alpha(z_n,z_n^*)}\to0.\]
Since $x\in D$, we can find a $\varepsilon$-neighborhood in $\mathbb{C}^+$ of $x$ on which $|m_{\calB}(z)|$ is bounded. While
\[ z_n - \gamma\int\frac{ud\calA(u)}{1+um(z_n)}  \to x,\]
it implies that for all sufficiently large $n$, $z_n- \gamma\int u(1+um(z_n))^{-1}d\calA(u)$ is in the $\varepsilon$-neighborhood of $x$.  
 Therefore, there exists a constant $C$ such that for all sufficiently large $n$ 
\[|m(z_n)| = \left|m_{\calB}\left( z_n- \gamma\int\frac{ud\calA(u)}{1+um(z_n)}\right)\right| <C,\]
which raise a contradiction. The proof is complete.
\end{proof}

\begin{lemma}\label{lemma:boundedness_integrals}
Consider any fixed $x\in \nonsig$. There exists $\varepsilon>0$ such that $\alpha(z,z^*)$ is bounded on $\frakN(x,\varepsilon)$.  
\end{lemma}
\begin{proof}
Suppose not. We can find a sequence of $z_n\in\mathbb{C}^+ \to x$ such that $\alpha(z_n, z_n^*)\to\infty$. Then, by Eq. \eqref{eq:ineq_core},
\[ \beta(z_n,z_n^*) = \gamma\int\frac{u^2d\calA(u)}{|1+um(z_n)|^2} < \frac{1}{\alpha(z_n,z_n^*)} \to 0.\]
It implies that the integrand $u^2/|1+um(z_n)|^2$ converges to $0$ $\calA$-almost everywhere as $n\to\infty$. Under the assumption that $\calA$ is not the Dirac measure at zero, we conclude that $|m(z_n)|\to\infty$, which contradicts with Lemma \ref{lemma:boundedness}.
\end{proof}

\begin{lemma}\label{lemma:boundedness_beta}
Consider any fixed $x\in \nonsig$. There exists $\varepsilon>0$ such that $\beta(z,z^*)$ is bounded on $\frakN(x,\varepsilon)$.
\end{lemma}
\begin{proof}
 Suppose not. We can find a sequence of $z_n\in\mathbb{C}^+\to0$ such that $\beta(z_n,z_n^*)\to\infty$. Using Eq. \eqref{eq:ineq_core}, we have $\alpha(z_n,z_n^*) < {1}/{\beta(z_n,z_n^*)}\to 0$.
It follows that $|m(z_n)|\leq (\alpha(z_n,z_n^*))^{1/2}\to0$. Here, the bound on $|m(z_n)|$ follows from taking modulus on both sides of Eq. \eqref{eq:main_eq} and  Cauchy-Schwarz. By the dominated convergence theorem
\[\beta(z_n,z_n^*) =\gamma \int \frac{u^2d\calA(u)}{|1+um(z_n)|^2} \to \gamma \int u^2d\calA(u) <\infty, \]
which raises a contradiction.
\end{proof}

\begin{lemma}\label{lemma:same_limit_converg}
Suppose that $\{z_n\}$ and $\{\tilde{z}_n\}$ are two sequences in $\mathbb{C}^+$ such that $z_n$ and $\tilde{z}_n$ both converging to $x\in \nonsig$ as $n\to\infty$. If $m_n = m(z_n)\to \underline{m}$ and $\tilde{m}_n =  m(\tilde{z}_n) \to \tilde{\underline{m}}$ as $n\to\infty$, then $\underline{m} = \tilde{\underline{m}}$. 
\end{lemma}
\begin{proof}
Suppose $\underline{m}\neq \tilde{\underline{m}}$. Using \eqref{eq:main_eq},
\begin{align*}
m_n - \tilde{m}_n &= \int\frac{d\calB(\tau)}{\tau - z_n + \gamma\displaystyle\int\frac{ud\calA(u)}{1+um_n}} - \int\frac{d\calB(\tau)}{\tau - \tilde{z}_n + \gamma\displaystyle\int\frac{ud\calA(u)}{1+u\tilde{m}_n}}\\
&= (z_n-\tilde{z}_n) \alpha(z_n,\tilde{z}_n) + (m_n-\tilde{m}_n) \alpha(z_n, \tilde{z}_n)\beta(z_n,\tilde{z}_n). 
\end{align*}
It follows that 
\begin{equation}\label{eq:diff_fraction}
m_n -\tilde{m}_n = \frac{ (z_n-\tilde{z}_n)\alpha(z_n, \tilde{z}_n)}{1 - \alpha(z_n,\tilde{z}_n)\beta(z_n,\tilde{z}_n)}.  
\end{equation}
By Cauchy-Schwarz and Lemma \ref{lemma:boundedness_integrals}, $0<|\alpha(z_n, \tilde{z}_n)| \leq \{\alpha(z_n, z_n^*) \alpha(\tilde{z}_n, \tilde{z}_n^*)\}^{1/2} <\infty$.
Therefore, $|(z_n- \tilde{z}_n)\alpha(z_n,\tilde{z}_n)|\to0$, as $n\to\infty$, since $z_n-\tilde{z}_n\to0$.  We shall show that 
$\liminf_{n\to\infty} \Big| 1-\alpha(z_n, \tilde{z}_n) \beta(z_n,\tilde{z}_n)\Big| >0$ and raise a contradiction. 
Note that 
\begin{align*}
&\Re\Big( \alpha(z_n,\tilde{z}_n) \beta(z_n,\tilde{z}_n)\Big) \\
&= \Re \int\int\frac{\gamma u^2d\calB(\tau)d\calA(u)}{\Big(\tau - z_n +\gamma\displaystyle\int\frac{wd\calA(w)}{1+wm_n}\Big)\Big(\tau - \tilde{z}_n +\gamma\displaystyle\int\frac{wd\calA(w)}{1+w\tilde{m}_n}\Big){(1+um_n)(1+u\tilde{m}_n)}}\\
& = \frac{1}{4}\gamma \int\int \left|\frac{u}{(1+um^*_n)\Big(\tau - z^*_n + \gamma\displaystyle\int\frac{wd\calA(w)}{1+wm^*_n}\Big)} + \frac{u}{(1+u\tilde{m}_n)\Big(\tau - \tilde{z}_n + \gamma\displaystyle\int\frac{wd\calA(w)}{1+w\tilde{m}_n}\Big)}   \right|^2d\calB(\tau)d\calA(u) \\
&\quad - \frac{1}{4}\gamma \int\int \left|\frac{u}{(1+um^*_n)\Big(\tau - z^*_n + \gamma\displaystyle\int\frac{wd\calA(w)}{1+wm^*_n}\Big)} - \frac{u}{(1+u\tilde{m}_n)\Big(\tau - \tilde{z}_n + \gamma\displaystyle\int\frac{wd\calA(w)}{1+w\tilde{m}_n}\Big)}   \right|^2d\calB(\tau)d\calA(u),\\
& = \frac{1}{4}(A_1 - A_2), \mbox{ say}.
\end{align*}
By Cauchy-Schwarz and $|a+b|^2 \leq 2|a|^2 + 2|b|^2$, 
\begin{align*}
A_1 \leq & 2 \gamma\int\frac{u^2d\calA(u)}{|1+um_n^*|^2} \int \frac{1}{\left|\tau - z^*_n + \gamma\displaystyle\int\frac{wd\calA(w)}{1+wm^*_n}\right|^2}d\calB(\tau)  + 2 \gamma\displaystyle\int\frac{u^2d\calA(u)}{|1+u\tilde{m}_n|^2} \int \frac{1}{\left|\tau - \tilde{z}_n + \gamma\displaystyle\int\frac{wd\calA(w)}{1+w\tilde{m}_n}\right|^2}d\calB(\tau) \\
=& 2 \alpha(z_n,z_n^*) \beta(z_n, z_n^*) + 2 \alpha(\tilde{z}_n, \tilde{z}_n^*) \beta(\tilde{z}_n, \tilde{z}_n^*)<4. 
\end{align*}
The last step is due to the bound $\alpha(z_n,z_n^*) \beta(z_n, z_n^*)<1$ by Eq. \eqref{eq:ineq_core}. 

If we can show $\liminf_{n\to\infty} A_2>0$, the proof will be complete, since
\begin{align*}
\liminf_{n\to\infty}\Big| 1-\alpha(z_n, \tilde{z}_n) \beta(z_n,\tilde{z}_n)\Big| \geq \liminf_{n\to\infty} 1 - \Re \alpha(z_n, \tilde{z}_n) \beta(z_n,\tilde{z}_n) > \liminf_{n\to\infty}\frac{1}{4} A_2>0.
\end{align*}
For succinctness, we present the proof of $\liminf_{n\to\infty} A_2>0$ in Appendix \ref{appendix:additional_proof_lemma_same_limit}. 
\end{proof}
Combining Lemma \ref{lemma:boundedness} and Lemma \ref{lemma:same_limit_converg}, we have shown that $m(z)$ is bounded in a neighborhood of $x\in\nonsig$ and any convergent subsequence converges to the same limit. The existence of $\underline{m}(x) = \lim_{z\in\mathbb{C}^+\to x}m(z)$ when $x\in D$ is proved. 

It remains to show the analyticity of the density function at any $x_0\in D$ such that $f(x_0) = \pi^{-1}\Im \underline{m}(x_0)>0$. 
Fix $x_0$ such that $\Im \underline{m}(x_0)>0$. Consider a neighborhood $\frakN$ of $(x_0, \underline{m}(x_0))$ in $\mathbb{C} \times \mathbb{C}^+$. We can select $\frakN$ to be sufficiently ``small'' so that $\frakN\subset \{(z, m): \Re z\in D, ~ ~|\Im z|\leq\kappa_1, \mbox{ and }\Im m \geq \kappa_2 \}$, where $\kappa_1$ and $\kappa_2$ are small positive constants. On $\frakN$, define 
\[ G(z, m)= \int\frac{d\calB(\tau)}{\tau -z + \gamma\displaystyle\int\frac{ud\calA(u)}{1+um}} -m .\]
Note that we can select $\kappa_1$ be sufficiently small so that $\Im z -\Im \gamma\int u(1+um)^{-1}d\calA(u)>0$ on $\frakN$. Therefore the integrand above is integrable. It is clear that $(z, m(z))$ is the solution to $G(z,m) = 0$ when restricted to $\mathbb{C}^+\times \mathbb{C}^+$. Moreover, at the origin
$G(x_0, \underline{m}(x_0)) = 0$, because by the dominated convergence theorem
\[\underline{m}(x_0) = \lim_{z\in\mathbb{C}^+\to x_0} m(z) = \lim_{z\in \mathbb{C}^+\to x_0} \int\frac{d\calB(\tau)}{\tau - z + \gamma\displaystyle\int\frac{ud\calA(u)}{1+um(z)} }  = \int \frac{d\calB(\tau)}{\tau - x_0 + \gamma\displaystyle\int\frac{ud\calA(u)}{1+u \underline{m}(x_0)}}.\]
Recall the holomorphic implicit function theorem. If we can show that 
\begin{equation}\label{eq:derivative_nonzero}
\frac{\partial G}{\partial m}(x_0,\underline{m}(x_0))\neq 0,
\end{equation}
the implicit function theorem suggests that there exists a neighborhood $\frakM$ of $x_0$, a neighborhood $\frakM'$ of $\underline{m}(x_0)$ so that $\frakM\times \frakM'\subset\frakN$,  and a holomorphic function $\underline{\tilde{m}}: \frakM\to\frakM'$ such that 
\[ \Big\{ (z,m)\in \frakM\times \frakM': G(z,m) =0 \Big\}  = \Big\{ (z, \underline{\tilde{m}}(z)): z\in\frakM\Big\}.\]
Since $\underline{\tilde{m}}(z)$ and $m(z)$ coincide on $z\in\frakM \cap \mathbb{C}^+$, we get that $\underline{\tilde{m}}(z)$ is a holomorphic extension of $m(z)$ on $\frakM$. It follows that $f(x_0) = \frac{1}{\pi}\Im \underline{\tilde{m}}(x_0)$ is analytic near $x_0$. Namely, $f(x) = \pi^{-1}\sum_{n=0}^\infty \Im a_n(x-x_0)^n$ for any real $x$ near $x_0$, where $a_n$ is a sequence of coefficients. 

It remains to show Eq. \eqref{eq:derivative_nonzero}. We call $\alpha(x_0,x_0) = \lim_{z\in\mathbb{C}^+ \to x_0}\alpha(z,z)$, $\beta(x_0,x_0) = \lim_{z\in\mathbb{C}^+ \to x_0}\beta(z,z)$, $A(x_0,x_0) = \lim_{z\in\mathbb{C}^+ \to x_0}\alpha(z,z^*)$, and $B(x_0,x_0) = \lim_{z\in\mathbb{C}^+ \to x_0}\beta(z,z^*)$. The existences of the limits follow from the dominated convergence theorem and the fact that $\Im\underline{m}(x_0)>0$. Using Eq. \eqref{eq:ineq_core}, 
\[ 1 -A(x_0,x_0)B(x_0,x_0) =  \lim_{z\to x_0} 1- \alpha(z, z^*)\beta(z,z^*) = \lim_{z\to x_0} \frac{\Im z}{\Im m(z)}\alpha(z,z^*) =0.\]
Recall that any integrable random variable $X$ satisfies $|\mathbb{E}X|\leq \mathbb{E} |X|$ and the equality holds if and only if $X = \theta |X|$ almost everywhere where $\theta$ is a modulus one constant. We therefore have
\begin{align*}
&|\alpha(x_0, x_0)| = \Big|\int \frac{d\calB(\tau)}{\left(\tau - x_0 + \gamma\displaystyle\int\frac{ud\calA(u)}{1+u\underline{m}(x_0)}\right)^2}\Big|
 < \int \frac{d\calB(\tau)}{\left|\tau - x_0 + \gamma\displaystyle\int\frac{ud\calA(u)}{1+u\underline{m}(x_0)}\right|^2} = A(x_0,x_0),\\[10pt]
&|\beta(x_0,x_0)| = \gamma \left|\int\frac{u^2d\calA(u)}{(1+u \underline{m}(x_0))^2}  \right| <  \gamma \int\frac{u^2d\calA(u)}{\left|1+u \underline{m}(x_0) \right|^2} = B(x_0,x_0).\\[10pt]
&\frac{\partial G}{\partial m}(x_0,\underline{m}(x_0)) = \int\frac{d\calB(\tau)}{\left( \tau -z +\gamma\displaystyle\int\frac{ud\calA(u)}{1+u\underline{m}(x_0)}\right)^2} \gamma\displaystyle\int\frac{u^2d\calA(u)}{(1+u\underline{m}(x_0))^2} -1 = \alpha(x_0,x_0)\beta(x_0,x_0)-1 \neq 0.
\end{align*}
The proof of Theorem \ref{thm:main} is complete. 

\subsection{Proof of Theorem \ref{thm:existence_of_masses}}\label{sec:existence_subsec:mass}

In this section, we show that the condition $\calB\{b\}-\gamma(1-\calA(\{0\}))>0$ is necessary for $F(\{b\})>0$. The sufficiency is shown in Appendix \ref{appendix:proof_theorem_existence_masses} for succinctness. Suppose that $b\in\mathbb{R}$ is any point such that $F(\{b\})>0$, due to Lemma \ref{lemma:point_mass_formula},
$\lim_{y\downarrow0}iy m(b+iy) = - F(\{b\})<0$. It implies that $\Im m(b+iy)$ diverges to $+\infty$ as $y\downarrow 0$. Therefore, 
\[ \left|\int\frac{ud\calA(u)}{1+um(b+iy) }\right| \leq \int_{u\neq0} \frac{d\calA(u)}{|\Im(m(b+iy))|} = \frac{1-\calA(\{0\})}{|\Im(m(b+iy))|} \to 0, \quad \mbox{as }y\downarrow 0. \]
Using Eq. \eqref{eq:main_eq}, 
\begin{align*}
iym(b+ iy) = \frac{iy\calB(\{b\})}{ -iy + \gamma\displaystyle\int\frac{ud\calA(u)}{1+um(b+iy)} } + \int_{\tau\neq b}\frac{iyd\calB(\tau)}{(\tau-b) -iy + \gamma\displaystyle\int\frac{ud\calA(u)}{1+um(b+iy)}}.
\end{align*}
For the second term, pointwise for any $\tau \neq b$,
\[ \frac{iy}{(\tau-b) -iy + \gamma\displaystyle\int\frac{ud\calA(u)}{1+um(b+iy)} } \longrightarrow 0, \quad \mbox{as }y\downarrow 0.\]
Moreover, uniformly for any $\tau$
\[ \left|\frac{iy}{(\tau-b) -iy + \gamma\displaystyle\int\frac{ud\calA(u)}{1+um(b+iy)} }\right| \leq \frac{|iy|}{ \left| \Im \Big(-iy + \gamma \displaystyle\int\frac{ud\calA(u)}{1+um(b+iy)}\Big)  \right| }\leq 1.\]
Here, we are using the fact that $\Im \int u(1+um(b+iy))^{-1} d\calA(u)\leq 0$. By the dominated convergence theorem,
\[ \lim_{y\downarrow0}\int_{\tau\neq b}\frac{iyd\calB(\tau)}{(\tau-b) -iy + \gamma\displaystyle\int\frac{ud\calA(u)}{1+um(b+iy)}} =0.\]
It follows that 
\begin{equation}\label{eq:expression_Fb}
 -F(\{b\}) = \lim_{y\downarrow0} iym(b+ iy) =\lim_{y\downarrow0} \frac{iy\calB(\{b\})}{ -iy + \gamma\displaystyle\int\frac{ud\calA(u)}{1+um(b+iy)}} =\lim_{y\downarrow0} \frac{\calB(\{b\})}{ -1 + \gamma\displaystyle\int\frac{ud\calA(u)}{iy+uiym(b+iy)}}. 
\end{equation}
Since $\lim_{y\downarrow 0}iy m(b+iy)<0$, we can find a constant $c>0$ such that $\Re(iym(b+iy))<-c$ for all sufficiently small $y$. Assume the condition is satisfied in the following. Therefore, for $u\neq 0$, 
\[  \left|\frac{u}{iy+uiym(b+iy)} \right| \leq \frac{|u|}{|u||\Re(iym(b+iy))|}\leq \frac{1}{c}.\]
Again by the dominated convergence theorem,
\begin{align*} 
\lim_{y\downarrow0}\gamma\int\frac{ud\calA(u)}{iy + uiym(b+iy)}& = \lim_{y\downarrow0} \gamma\int_{u\neq 0}\frac{ud\calA(u)}{iy+uiym(b+iy)} =  \gamma\int_{u\neq 0}\lim_{y\downarrow0}\frac{ud\calA(u)}{iy+uiym(b+iy)}\\[10pt]
&  =\frac{\gamma (1-\calA(\{0\}))}{\lim_{y\downarrow0}iym(b+iy)}  =\frac{\gamma (1-\calA(\{0\}))}{-F(\{b\})}  .
\end{align*}
Plugging into Eq. \eqref{eq:expression_Fb}, we have 
\[ -F(\{b \}) = \frac{\calB(\{b\})}{-1 - \gamma(1-\calA(\{0\}))/F(\{b\}) }. \]
That is, $F(\{b\}) = \calB(\{b\}) - \gamma(1-\calA(\{0\}))>0$. 



\section{Support}\label{sec:support}


For a distribution $G$, recall that its support is \(\supp(G) = \{x \in \mathbb{R} : G(U) > 0 \text{ for every open set } U \text{ containing } x\}\). While Theorem \ref{thm:main} indicates that the density of \(F\) on \(D\) exists and Theorem \ref{thm:existence_of_masses} explains the mass of $F$ on any point, the results shed little light on the support of $F$. In this section, we present the characterization of the support of \(F\) given \(\calA\) and \(\calB\).

It is more convenient to consider the complement of the support. Notice that if \(x \in \supp^\complement(F)\),
\[
\lim_{z \in \mathbb{C}^+ \to x} m(z) = \lim_{z \in \mathbb{C}^+ \to x}\int \frac{dF(\tau)}{\tau - z} =\int \frac{dF(\tau)}{\tau - x}= \underline{m}(x) \text{ exists and is real}.
\]
Moreover, \(\underline{m}(x)\) is analytic and monotonically increasing  on any connected component in \(\supp^\complement(F)\), since 
\[
\underline{m}'(x) = \int \frac{dF(\tau)}{(\tau - x)^2} > 0.
\]
These results actually hold true for the Stieltjes transform of any probability measure.



\subsection{Properties of the Stieltjes transform outside the support of $F$}\label{sec:support_subsec:lemmas}
We first describe the range of $\underline{m}(x)$ when $x\in \supp^\complement(F)$. Define
\begin{equation}\label{eq:def_EcalA}
E(\calA) = \begin{cases}
\{0\}\cup\Big\{m\neq 0: -m^{-1} \notin \supp(\calA)\Big\}, & \mbox{if }\supp(\calA)\mbox{ is compact}; \\ 
\Big\{ m\neq 0: -m^{-1} \notin \supp(\calA)\Big\}, & \mbox{otherwise}.
 \end{cases}
\end{equation}
The integrand of $\int u(1+um)^{-1}d\calA(u)$ is well-defined and integrable for $m\in E(\calA)$. 
\begin{lemma}\label{lemma:existence_extension_calA}
Consider any $x_0\in\supp^{\complement}(F)$. We have $\underline{m}(x_0)\in E(\calA)$. 
Moreover, define $h(x_0) = x_0 - \gamma \int{u(1+u\underline{m}(x_0))^{-1}d\calA(u)}$.
Then, $h(x_0)\notin \supp(\calB)$. 
\end{lemma}

\begin{proof}
We first show that if $x_0\in\supp^\complement(F)$ and $\underline{m}(x_0)\neq 0$, we have that $-\underline{m}^{-1}(x_0) \notin \supp(\calA)$. 
Since $\supp^{\complement}(F)$ is open and $\underline{m}(x)$ is continuous, there exists an open interval containing $x_0$, denoted by $(x_1, x_2)$, such that $(x_1, x_2)\subset\supp^\complement(F)$ and 
\[   \underline{m}(x) =  \lim_{y\downarrow0} m(x+iy) \in \mathbb{R}\setminus\{0\}, \quad x\in (x_1, x_2).\] 
We show that there exists $y_0>0$ such that $\beta(z, z^*)$ is bounded on $\{z = x+iy: x\in (x_1,x_2),~0< y<y_0 \}$.  Suppose not. We can find a sequence of $z_n\to x\in(x_1,x_2)$ such that $\beta(z_n,z_n^*)\to \infty$. 
Using Eq. \eqref{eq:ineq_core}, $\alpha(z_n,z_n^*)\to0$. It follows that $|m(z_n)|\leq \alpha^{1/2}(z_n,z_n^*) \to 0$. It is a contradiction with the assumption $\underline{m}(x)\neq 0$. 
By the Fatou's lemma, 
\[ \int \frac{|u|}{|1+u\underline{m}(x)|}d\calA(u) \leq\liminf_{n\to\infty} \left(\int \frac{|u|^2}{|1+u m(z_n)|^2}d\calA(u)\right)^{1/2}<\infty.\]
By the dominated convergence theorem, if $z\in\mathbb{C}^+\to x\in(x_1,x_2)$, $\int u(1+um(z))^{-1}d\calA(u) \to \int u(1+u\underline{m}(x))^{-1}d\calA(u)$. We have then
\[ m_{\calA}(-1/m(z))  = m(z) + m^2(z) \int\frac{ud\calA(u)}{ 1+ um(z) } \to \underline{m}(x) + \underline{m}^2(x) \int\frac{ud\calA(u)}{ 1+ u\underline{m}(x)}.\]
It indicates that $m_{\calA}$ is real in a neighborhood of $-\underline{m}^{-1}(x_0)$, which implies that $-1/\underline{m}(x_0)\notin\supp(\calA)$. 

Next, we show that if there exists $x_0\in\supp^{\complement}(F)$ and $\underline{m}(x_0)=0$, then $\supp(\calA)$ is compact. Since $\supp^{\complement}(F)$ is open, we can find  $(x_1, x_2)$ containing $x_0$ such that $(x_1, x_2)\subset\supp^\complement(F)$. While $\underline{m}'(x) >0$ on $(x_1, x_2)$, we conclude that as $x\uparrow x_0$, $\underline{m}(x) \uparrow 0$ and as $x\downarrow x_0$, $\underline{m}(x) \downarrow 0$. Using the previous arguments, for $x\neq x_0$ and $x\in(x_1, x_2)$, $-1/\underline{m}(x) \notin \supp(\calA)$. Since as $x\uparrow x_0$, $-1/\underline{m}(x) \uparrow+\infty$ and as $x\downarrow x_0$,  $-1/\underline{m}(x) \downarrow -\infty$, we conclude that there exists $a_1$ and $a_2$ such that $(-\infty, a_1) \subset \supp^\complement(\calA)$ and $(a_2, +\infty) \subset \supp^\complement(\calA)$. That is, $\supp(\calA)$ is compact. Combining all these arguments, we have that $\underline{m}(x_0) \in E(\calA)$.

Lastly, we show that $h(x_0)\notin \supp(\calB)$. Following the previous statement, there exists $y_0>0$ such that $\alpha(z,z^*)$ is bounded on $\{z = x+iy: x\in (x_1,x_2), ~0< y<y_0 \}$. Suppose not. We can find a sequence $z_n\to x\in (x_1,x_2)$ such that $\alpha(z_n,z_n^*)\to\infty$. Using Eq. \eqref{eq:ineq_core}, $\beta(z_n,z_n^*)\to 0$. It follows that 
$u^2/|1+um(z_n)|^2 \to 0$ $\calA$-almost everywhere, as $z_n\to x$. Given that $\calA$ is not the Dirac measure at $0$, it implies $|m(z_n)|\to\infty$, which contradicts with the assumption that $\underline{m}(x)$ exists and is real. Using the dominated convergence theorem, for any $x\in (x_1, x_2)$,
\[\underline{m}(x) = \lim_{z\in\mathbb{C}^+\to x} m(z) = \lim_{z\in\mathbb{C}^+\to x}\int\frac{d\calB(\tau)}{\tau - z + \gamma \displaystyle\int\frac{ud\calA(u)}{1+um(z)}} = m_{\calB}\left(x -\gamma \displaystyle\int\frac{ud\calA(u)}{1+u\underline{m}(x)}  \right).\]
It indicates that $m_{\calB}$ is real at $x - \gamma \int u(1+u\underline{m}(x))^{-1}d\calA(u)$. Moreover, on $(x_1, x_2)$,  $x- \gamma \int u(1+u\underline{m}(x))^{-1} d\calA(u)$ is monotonically increasing as 
\[ 1 +  \underline{m}'(x)\gamma \int \frac{u^2d\calA(u)}{(1+u\underline{m}(x))^2} >0.\] 
It implies that $m_{\calB}$ is real on a neighborhood of $h(x_0)$. Consequently, $h(x_0)\notin\supp(\calB)$. 
\end{proof}

\subsection{Determination of the support of $F$}

Consider all possible \((x, m) \in \mathbb{R}^2\) such that 
\begin{equation}\label{eq:extension_main_eq_to_real_line}
\begin{cases}
&m \in E(\calA),\\[10pt]
&x - \gamma \displaystyle\int \frac{u \, d\calA(u)}{1 + u m} \in \supp^\complement(\calB),\\[10pt]
&m = \displaystyle\int \frac{d\calB(\tau)}{\tau - x + \gamma \displaystyle\int \frac{u \, d\calA(u)}{1 + u m}}.
\end{cases}
\end{equation}
Lemma \ref{lemma:existence_extension_calA} indicates that when \(x_0 \in \supp^\complement(F)\), Eq. \eqref{eq:extension_main_eq_to_real_line} is satisfied at \((x_0, \underline{m}(x_0))\). However, it is not hard to discover that for  general \(x_0 \in \mathbb{R}\), there can be multiple \(m\) values satisfying Eq. \eqref{eq:extension_main_eq_to_real_line}. For example, when \(\calB\) and \(\calA\) are both Dirac measures, the last equation becomes a polynomial in \(m\), and multiple real roots may exist. Readers should not confuse Eq. \eqref{eq:extension_main_eq_to_real_line} with Eq. \eqref{eq:main_eq}, as Eq. \eqref{eq:main_eq} maps \(\mathbb{C}^+\) to \(\mathbb{C}^+\), and its solution \(m(z)\) on $\mathbb{C}^+$ is unique for each \(z\). Conversely, for some \(m \in \mathbb{R}\), there are possibly multiple \(x \in \mathbb{R}\) satisfying Eq. \eqref{eq:extension_main_eq_to_real_line}. Moreover, the solution to the equation exists for some \(x\) beyond \(\supp^\complement(F)\). 

For convenience, we reformulate Eq. \eqref{eq:extension_main_eq_to_real_line} as follows. Define
\[ 
\calH = \Big\{ h \in \supp^\complement(\calB) : m_{\calB}(h) = \int (\tau - h)^{-1} d\calB(\tau) \in E(\calA) \Big\}.
\]
Note that $\calH$ is open.  For all \(h \in \calH\), let 
\begin{equation}\label{eq:system_h_m_x} 
\begin{cases}
m_{\calB}(h) = \displaystyle\int \frac{d\calB(\tau)}{\tau - h},\\[10pt]
x = x(h) = h + \gamma \displaystyle\int \frac{u \, d\calA(u)}{1 + u m_{\calB}(h)}.
\end{cases}
\end{equation}
The roots to the two systems are equivalent in the sense that 
\[
\Big\{\text{all possible solutions } (x, m) \text{ to Eq. \eqref{eq:extension_main_eq_to_real_line}} \Big\} = \Big\{(x(h), m_{\calB}(h)) : h \in \calH \Big\}.
\]
The advantage of Eq. \eqref{eq:system_h_m_x} is that for each \(h \in \calH\), there is a unique \(x\) satisfying Eq. \eqref{eq:system_h_m_x}, and the expression for \(x(h)\) is simple.

To characterize the support of $F$, we only need to determine for what values of $h$, $x(h)\in \supp^\complement(F)$. 
The following result shows that it is fully characterized by the sign of the derivative of $x$ with respect to $h$. 
\begin{theorem}\label{thm:increasing_outside_necessary_sufficient}
Suppose that $x_0\in \supp^\complement(F)$. Then, there is a unique $h_0\in\calH$ such that $(x(h_0), m_{\calB}(h_0)) = (x_0, \underline{m}(x_0))$, and 
\begin{equation}\label{eq:derivative_x_to_h_positive} 
x'(h_0) = 1 - \gamma\int\frac{u^2d\calA(u)}{(1+um_{\calB}(h_0))^2}\int\frac{d\calB(\tau)}{(\tau-h_0)^2} >0. 
\end{equation}
Reversely, suppose $h_0\in\calH$ and 
  Eq. \eqref{eq:derivative_x_to_h_positive} is satisfied. 
Then, $x(h_0) \in \supp^\complement(F)$  and $\underline{m}(x(h_0)) = m_{\calB}(h_0)$.   
\end{theorem}

\begin{remark}\label{remark:equivalence_to_Silverstein_Dozier_Coull}
Although \(h \mapsto m_{\calB}(h)\) is not one-to-one globally for general \(\calB\), the inverse mapping can be defined locally in a neighborhood of \(m_0 = m_{\calB}(h_0)\) since \(m'_{\calB}(h_0) > 0\) for any \(h_0 \in \supp^\complement(\calB)\). We have then locally
\[
\left.\frac{dx}{dm}\right|_{m = m_0} = \left.\frac{dx}{dh}\right|_{h = h_0} \left.\frac{dh}{dm}\right|_{m = m_0}.
\]
In \cite{silverstein1995analysis}, \cite{dozier2007analysis}, and \cite{couillet2014analysis}, the sign of the derivative of the inverse of the Stieltjes transform, that is, \(dx/dm\) under our settings, is used to determine the support of the limiting spectral distribution of a random matrix. For details, see Theorem 4.1 and Theorem 4.2 of \cite{silverstein1995analysis}, Theorem 3.2 and Theorem 3.3 of \cite{dozier2007analysis}, Section 3.2 and Section 3.3 of \cite{couillet2014analysis}. Since \(dh/dm > 0\), \(dx/dh\) has the same sign as \(dx/dm\). It indicates that our method follows the same idea as their methods but is represented in a different form. Our presentation has the advantage that the inverse function \(x(h)\) exists and exhibits an explicit form, while the inverse Stieltjes transform is not one-to-one and can only be defined locally.

\end{remark}

\begin{proof}[Proof of Theorem \ref{thm:increasing_outside_necessary_sufficient}]
Suppose that $x_0\in \supp^\complement(F)$. By Lemma \ref{lemma:existence_extension_calA}, Eq. \eqref{eq:main_eq} can be extended to $z= x_0$ as 
\begin{equation}\label{eq:proof_theorem_increasing_outside_1}
\underline{m}(x_0)  = \int \frac{d\calB(\tau)}{\tau  -x_0 + \gamma\displaystyle\int\frac{ud\calA(u)}{1+u\underline{m}(x_0)}}.
\end{equation}
Let $h_0 = x_0 - \gamma\int{u(1+u\underline{m}(x_0))^{-1}d\calA(u)}$. Then, $h_0\in\calH$ and  $(x(h_0), m_\calB(h_0)) = (x_0, \underline{m}(h_0))$. To show the uniqueness, suppose there exists  $h_1\in\calH \neq h_0$ such that $(x(h_1), m_\calB(h_1)) = (x_0, \underline{m}(x_0))$. Then,
\[ h_1 = x(h_1) - \gamma\int\frac{ud\calA(u)}{1+um_{\calB}(h_1)} = x_0 -\gamma\int\frac{ud\calA(u)}{1+u\underline{m}(x_0)} = x(h_0)-\gamma\int\frac{ud\calA(u)}{1+um_{\calB}(h_0)} = h_0,\]
which raises a contradiction. To show $x'(h_0)$ is positive, note that since $\underline{m}(x_0)$ is the Stieltjes transform of $F$ at $x_0\in\supp^{\complement}(F)$, $\underline{m}'(x_0)>0$. Differentiating both sides of Eq. \eqref{eq:proof_theorem_increasing_outside_1},
\[\underline{m}'(x_0) = \int \frac{d\calB(\tau)}{\left(\tau - x_0 +\gamma \displaystyle\int\frac{ud\calA(u)}{1+u\underline{m}(x_0)}\right)^2} \left(1 + \gamma\int \frac{u^2d\calA(u)}{(1+u\underline{m}(x_0))^2}m'(x_0)\right).\]
It follows that
\[ x'(h_0) = 1 - \gamma\int\frac{u^2d\calA(u)}{(1+um_{\calB}(h_0))^2}\int\frac{d\calB(\tau)}{(\tau-h_0)^2} = \frac{1}{\underline{m}'(x_0)}\int\frac{d\calB(\tau)}{\left(\tau -x_0 +\gamma\displaystyle \int \frac{ud\calA(u)}{1+u\underline{m}(x_0)}\right)^2}>0.\]

Reversely, suppose that $h_0 \in\calH$ and $x'(h_0)>0$. We show that $x_0 = x(h_0)\in\supp^\complement(F)$. It suffices to show that $\lim_{z\in\mathbb{C}^+\to x_0}m(z)$ exists and is real. Consider a neighborhood of $h_0$, say $(h_1, h_2)$, such that $(h_1, h_2) \subset \calH$ and $x'(h)>x'(h_0)/2$ for all $h\in(h_1,h_2)$. Such a neighborhood exists since $\calH$ is open and $x'(h)$ is continuous on $\calH$.

We first extend Eq. \eqref{eq:system_h_m_x} to the complex domain. For $|y|< y_0$, call  
\[\mathcal{O} = \{ (h, y): h\in(h_1,h_2), |y|<y_0\}.\] 
For $(h,y)\in \mathcal{O}$, consider
\begin{align*}
&\tilde{m}(h,y) = \int \frac{d\calB(\tau)}{\tau - h - iy},\\
&z(h,y) = h +iy + \gamma \int \frac{ud\calA(u)}{1+u\tilde{m}(h,y)}.
\end{align*}
Clearly, $(z(h,y), \tilde{m}(h,y))$ is holomorphic on $\mathcal{O}$. If $y= 0$, for any $h\in (h_1, h_2)$.
\[1 - \gamma \int\frac{u^2 d\calA(u)}{ |1 + u\tilde{m}(h,0)|^2}  \int\frac{d\calB(\tau)}{|\tau -h |^2} = x'(h)>x'(h_0)/2>0.\]
We can select $y_0$ to be sufficiently small such that on $\mathcal{O}$,
\[  1 - \gamma \int\frac{u^2 d\calA(u)}{ |1 + u\tilde{m}(h,y )|^2}  \int\frac{d\calB(\tau)}{|\tau - h -iy |^2} >0,\]
because the function is clearly continuous in $(h,y)$ near $(h_0, 0)$. In the following, assume this condition on $y_0$ is satisfied. 
Then, $\Im z(h,y)$ has the same sign as $y$ since 
\[ \Im z(h,y) = y \left(1 - \gamma \int\frac{u^2 d\calA(u)}{ |1 + u\tilde{m}(h,y)|^2}  \int\frac{d\calB(\tau)}{|\tau -h - iy|^2}\right).\]
Moreover, $\Im \tilde{m}(h,y)  = y \int |\tau -h -iy|^{-2}\calB(\tau)$ also has the same sign as $y$. Then,
\[\Big\{ \Big(z(h,y), \tilde{m}(h,y)\Big):~ h\in (h_1, h_2), ~0< y<y_0\Big\} \subset \mathbb{C}^+ \times \mathbb{C}^+.\]
If we restrict $y$ to be positive, 
\[ \Big\{z(h,y):  (h,y)\in \mathcal{O}^+\Big\} \subset \mathbb{C}^+,\]
where $\mathcal{O}^+ =\{ (h, y): h\in(h_1,h_2), 0<y<y_0\}$.

Next, we discuss the image set of $z(h,y)$ when  $(h,y)\in \mathcal{O}^+$. By the open mapping theorem, since $z(h,y)$ is a non-constant holomorphic function, it maps the open set $\mathcal{O}$ into an open set in $\mathbb{C}$. Clearly, $x_0 = z(h_0, 0)$ is in the image set.  Therefore, we can find an open ball $\frakN$ in $\mathbb{C}$ of $x_0$ such that $\frakN \subset \Big\{ z(h,y): (h,y)\in\mathcal{O}\Big\}$ and the half open ball $\frakN \cap \mathbb{C}^+ \subset \Big\{ z(h,y): (h,y)\in\mathcal{O}^+\Big\}$.

Lastly, we show that $\Big\{(z(h,y), \tilde{m}(h,y)): (h,y)\in\mathcal{O}\Big\}$ is a holomorphic extension of the Stieltjes transform $(z, m(z))$ on $z\in \Big\{z(h,y): (h,y)\in\mathcal{O}^+ \Big\} $. It is indeed quite straightforward since by the definition of $z(h,y)$ and $\tilde{m}(h,y)$, we have  that $\{ (z(h,y), \tilde{m}(h,y)): (h,y)\in\mathcal{O}\}$ satifies Eq. \eqref{eq:main_eq}, namely 
\[\tilde{m}(h,y) = \int\frac{d\calB(\tau)}{\tau - z(h,y) + \gamma\displaystyle\int\frac{ud\calA(u)}{1+u\tilde{m}(h,y)}}.\]
Recall that the solution in $\mathbb{C}^+$ to Eq. \eqref{eq:main_eq} is unique for any $z\in\mathbb{C}^+$. 
We must have $(z(h,y), \tilde{m}(h,y))$, $(h,y)\in \mathcal{O}^+$ match with $(z(h,y), m(z(h,y))$,  $(h,y)\in\mathcal{O}^+$. Our conclusion follows.

As the half open ball $\frakN\cap \mathbb{C}^+$  is contained in the image set of $\{z(h,y), (h,y)\in\mathcal{O}^+\} $, $z(h,y)$ can approach $x_0$ from all directions in $\mathbb{C}^+$. We find that
\[ \lim_{z\in\mathbb{C}^+ \to x_0} m(z)  = \tilde{m}(h_0, 0)  =\int \frac{d\calB(\tau)}{\tau - h_0}.\]
As $\tilde{m}(h_0, 0)$ is real, it follows that $x_0 \in \supp^\complement(F)$. 
\end{proof}

\subsection{Practical procedure for determining the support}
\begin{figure}[htbp]
\centering 
\includegraphics[width = \textwidth]{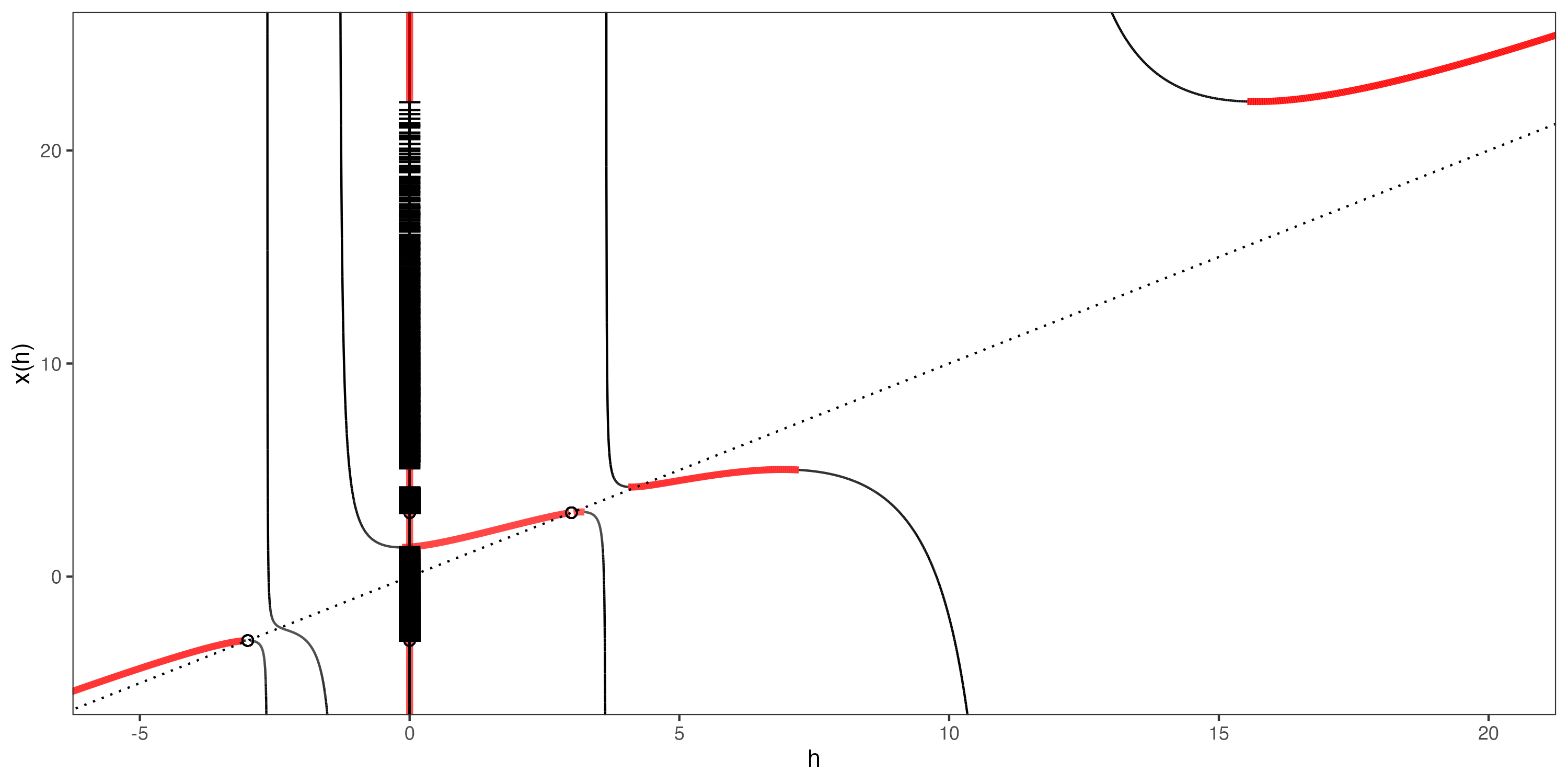}
\caption{$x(h)$ versus $h$. Settings: $\calA = 0.2 \delta_0 + 0.4 \delta_1 + 0.4 \delta_{10}$, $\calB = 0.4\delta_{-3} + 0.6 \delta_{3}$, and $\gamma =0.5$. In red thick line, positions for which $x'(h)>0$. On the vertical axis, the red segments are $\supp^\complement(F)$ (the projection of the red thick lines to the axis) and what left is $\supp(F)$. In black dashes on the axis, empirical eigenvalue points for $p=1000$. The two hollow dots are $(-3,-3)$ and $(3,3)$. Note $x=3$ is an isolated point of $\supp(F)$ and $x=-3\in\partial\supp(F)$ but $F(\{-3\})=0$. The dotted slope is $x=h$. }
\label{fig:fig1}
\end{figure}

\begin{figure}[htbp]
\centering 
\includegraphics[width = \textwidth]{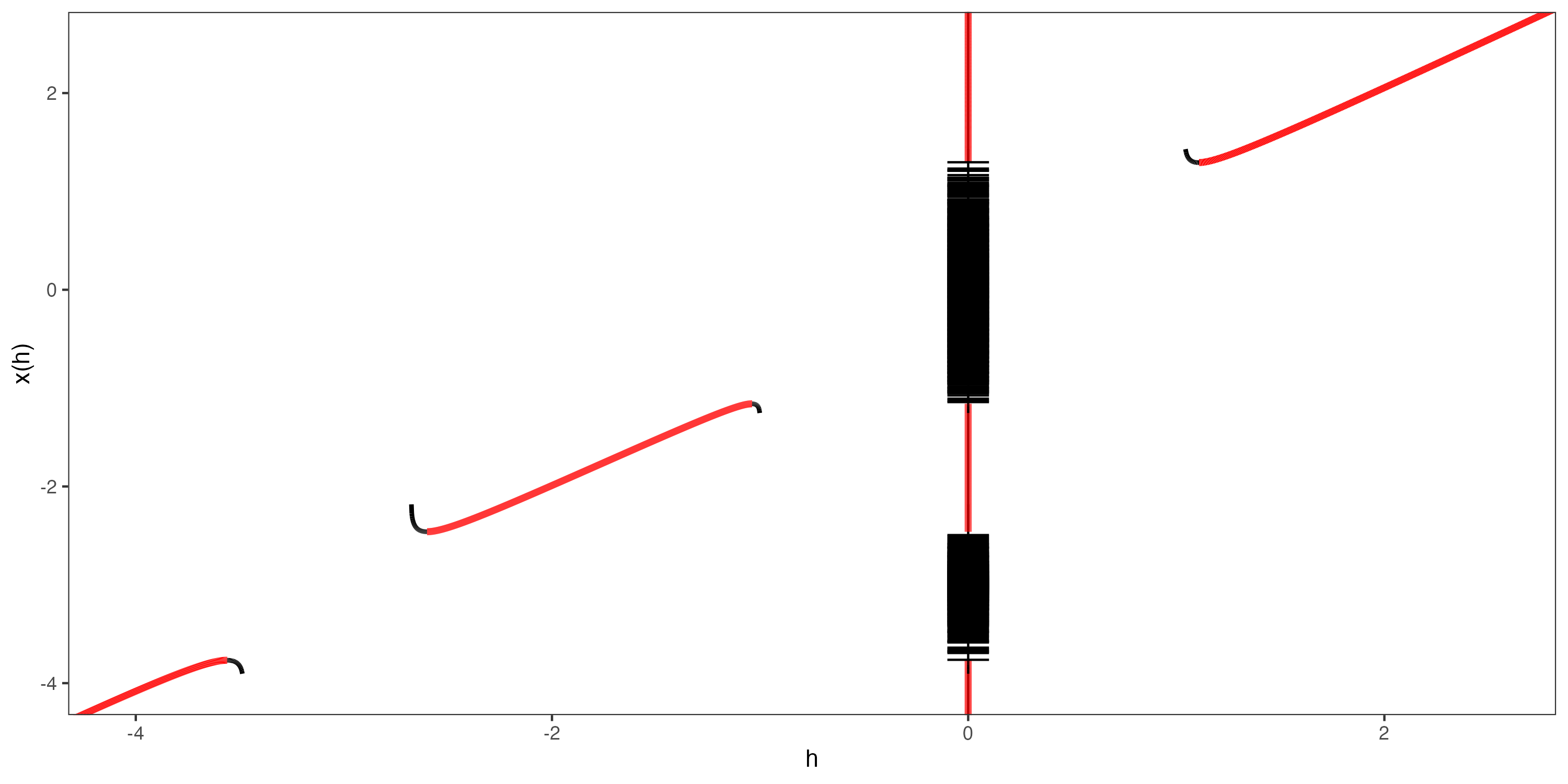}
\caption{$x(h)$ versus $h$. Settings: $\calA = SC_1$, $\calB = 0.4\delta_{-3} +0.6SC_1$ and $\gamma=0.5$, where $SC_1$ is the semicircular law with radius $1$. In red thick line, positions for which $x'(h)>0$. On the vertical axis, the red segments are $\supp^\complement(F)$ (the projection of the red thick lines to the axis) and what left is $\supp(F)$. In black dashes on the axis, empirical eigenvalue positions for $p=400$. }
\label{fig:fig2}
\end{figure}

Following Theorem \ref{thm:increasing_outside_necessary_sufficient}, we can use the graph of $x$ against $h$ as a tool to determine the support of $F$. The procedure is summarized in Algorithm \ref{algo:determination_support_F}. 
Two examples are given in Figure \ref{fig:fig1} and Figure \ref{fig:fig2}. 

\begin{algorithm}
\caption{Determining the Support of Distribution $F$}\label{algo:determination_support_F}
\begin{algorithmic}[1]
\REQUIRE Probability distributions $\mathcal{A}$, $\mathcal{B}$, and a constant $\gamma>0$.
\ENSURE The support of $F$ denoted by $\supp(F)$. 
\STATE If $\supp(\calB) = \mathbb{R}$, then return $\supp(F) = \mathbb{R}$. If not, perform the following steps. 
\STATE  Discard $\supp(\calB)$ on the horizontal axis. Let $h$ run through $\text{supp}^\complement(\mathcal{B})$ and calculate the Stieltjes transform:
\[
m_{\mathcal{B}}(h) = \int \frac{1}{\tau - h} d\mathcal{B}(\tau).
\]
Further, if the support of $\calA$ is compact, discard on the horizontal axis all $h$ such that $m_{\calB}(h)\neq 0$ and $-m^{-1}_{\mathcal{B}}(h) \in \text{supp}(\mathcal{A})$; otherwise, discard all $h$ such that $m_{\calB}(h)=0$ or   $-m^{-1}_{\mathcal{B}}(h) \in \text{supp}(\mathcal{A})$.
 What remains on the horizontal axis is $\calH$. If $\calH= \emptyset$, return $\supp(F) = \mathbb{R}$. If not, perform the following steps. 
\STATE  Plot $x$ against $h$ on $\calH$ where
\[
x(h) = h + \gamma \int \frac{u}{1 + u m_{\mathcal{B}}(h)} d\mathcal{A}(u).
\]
\STATE  Find all points $h$ (if any) at which $x(h)$ is increasing. For all these $h$, remove $x(h)$ from the vertical axis. What remains on the vertical axis is $\supp(F)$. 
\end{algorithmic}
\end{algorithm}

Each of the following propositions in this section sheds some light on the qualitative behavior of the graph of $x$ versus $h$.
\begin{proposition}
\label{prop:x_h_increasing}
Suppose that $h_1\in\calH$ and $h_2\in\calH$ are such that $h_1<h_2$, $x'(h_1)>0$, $x'(h_2)>0$. Then, $x(h_1)<x(h_2)$. 
\end{proposition}
\begin{proof}
Using Eq. \eqref{eq:system_h_m_x},
\begin{align*}
x(h_1) - x(h_2) &= h_1 - h_2 + \gamma\int \left( \frac{u}{1+um_\calB(h_1)} - \frac{u}{1+um_\calB(h_2)}\right)d\calA(u)\\
& = (h_1-h_2)\left( 1-\gamma \int\frac{u^2d\calA(u)}{(1+um_{\calB}(h_1))(1+um_{\calB}(h_2))}\int\frac{d\calB(\tau)}{(\tau-h_1)(\tau-h_2)} \right).
\end{align*}
Using Cauchy-Schwarz, 
\[\left|\gamma \int\frac{u^2d\calA(u)}{(1+um_{\calB}(h_1))(1+um_{\calB}(h_2))}\int\frac{d\calB(\tau)}{(\tau-h_1)(\tau-h_2)}\right| \leq (1- x'(h_1))^{1/2}(1-x'(h_2))^{1/2} <1. \]
Therefore, $x(h_1) - x(h_2)$ has the same sign as $h_1 - h_2$. 
\end{proof}
\begin{proposition}
\label{prop:no_decreasing_in_middle}
Suppose that $[h_1, h_2] \subset\calH$. If $x'(h_1) >0$ and $x'(h_2)>0$, then $x'(h)>0$ for all $h\in [h_1,h_2]$.
\end{proposition}
\begin{proof}
Without loss of generality, take $h_1 <h_2$. By Proposition \ref{prop:x_h_increasing}, $x(h_1) <x(h_2)$. Since $x(h)$ is differentiable on $[h_1, h_2]$, for any $x_0\in[x(h_1), x(h_2)]$, we can find $h_0\in[h_1, h_2]$ such that $x_0 = x(h_0)$ and $x'(h_0)>0$. Therefore, $[x(h_1), x(h_2)] \subset \supp^\complement(F)$ by Theorem \ref{thm:increasing_outside_necessary_sufficient}. We get that $\underline{m}'(x) >0$ on $[x(h_1), x(h_2)]$. Recall that $x'(h)$ has the same sign as $dx/dm$ as explained in Remark \ref{remark:equivalence_to_Silverstein_Dozier_Coull}. It follows that $x'(h)>0$, for all $h\in[h_1,h_2]$.
\end{proof}
\begin{proposition}\label{prop:x_as_h_approach_infty}
If $\calA$ and $\calB$ are both compact, we can find $h_1$ and $h_2$ such that  $(-\infty, h_1) \cup (h_2, +\infty) \subset \calH$. Moreover, as $h\to \pm\infty$, $x(h)/h \to 1$ and $x'(h)\to 1$. 
\end{proposition}
\begin{proof}
If $\calB$ is compact, as $h\to +\infty$, $m_{\calB}(h) \uparrow 0$ and $-1/m_{\calB}(h) \uparrow \infty$. Therefore, we can find $h_2$ such that for $h> h_2$, $h\notin\supp(\calB)$ and  $-1/m_\calB(h) \notin \supp(\calA)$. As $h\to+\infty$,  since $hm_{\calB}(h)\to -1$ and $h^2m'_{\calB}(h)\to 1$,
\[ x(h)/h  =  1+ \int \frac{u d\calA(u)}{h + u hm_{\calB}(h)}\to1.\]
\[ x'(h) = 1 - \int\frac{h^2 d\calB(\tau)}{(\tau - h)^2} \gamma \int\frac{u^2d\calA(u)}{(h+ uhm_\calB(h))^2 } \to 1.\]
The case where $h\to-\infty$ can be treated similarly. 
\end{proof}
\begin{proposition}\label{prop:when_h_approaches_discontinuous_point_of_B}
Suppose that $h_0$ is an isolated point of $\supp(\calB)$ and there exists $\varepsilon>0$ such that $(h_0-\varepsilon, h_0) \cup (h_0, h_0 +\varepsilon) \in\calH$, we have $\lim_{h\in\calH\to h_0} x(h)= h_0$ and $\lim_{h\in\calH\to h_0}x'(h) = 1-  \gamma(1-\calA(\{0\}))/\calB(\{h_0\})$.
\end{proposition}
\begin{proof}
As $h\in\calH \to h_0$, $ |m_{\calB}(h)| \to \infty$, $(h_0- h)m_{\calB}(h) \to \calB(\{h_0\})$ and $(h_0- h)^2\int (\tau-h)^{-2} d\calB(\tau) \to \calB(\{h_0\})$. Therefore, by the dominated convergence theorem, as $h\to h_0$, $x(h) = h + \gamma \int u(1+um_{\calB}(h))^{-1}d\calA(u) \to h_0$, and 
\begin{align*}
x'(h)& = 1- \gamma \int\frac{u^2d\calA(u)}{(1+um_{\calB}(h))^2} \int \frac{d\calB(\tau)}{(\tau -h)^2}  = 1- \gamma\int\frac{u^2d\calA(u)}{( (h_0-h) +u (h_0-h)m_{\calB}(h))^2} \int \frac{(h_0-h)^2 d\calB(\tau)}{(\tau -h)^2}\\
&\to 1- \gamma (1-\calA(\{0\}))\frac{\calB(\{h_0\})}{(\calB(\{h_0\}))^2 } = 1 - \frac{\gamma (1-\calA(\{0\}))}{\calB(\{h_0\})}.
\end{align*}

\end{proof}

\section{Behavior near a boundary point}\label{sec:boundary_behavior}

We now study the behavior of \(F\) near a boundary point \(x_0\) of \(\supp(F)\). Following the assumption in \cite{silverstein1995analysis}, \cite{dozier2007analysis}, and \cite{couillet2014analysis}, we restrict to the case where \(x_0 = x(h_0)\) for some \(h_0 \in \calH\). While \(\calH\) is open, in that case, we clearly have \(x(h)\) being analytic around \(h_0\) and \(x'(h_0) = 0\). Theorem \ref{thm:increasing_outside_necessary_sufficient} indicates that when \(x_0\) is a left endpoint, \(x(h_0)\) is a local supremum. When \(x_0\) is a right endpoint, \(x(h_0)\) is a local infimum.

This assumption is not satisfied for some choices of \(\calA\) and \(\calB\). For instance, we can have \(x(h) \uparrow x_0\) as \(h \in \calH \uparrow h_0\), for some \(h_0 \in \partial \calH\). The two points \(x_0 = -3\) and \(x_0 = 3\) in Figure \ref{fig:fig1} are such boundary points of \(\supp(F)\). Indeed, \(x_0 = 3\) is an isolated point of \(\supp(F)\), and \(x_0 = -3\) is a left endpoint of an open subset of \(\supp(F)\). We note, however, that the assumption holds when both \(\calA\) and \(\calB\) are discrete and \(x_0\) is not a discontinuous point of \(\calB\) such that  $\calB(\{x_0\}) =\gamma(1-\calA(\{0\}))$.

Recall the definition of $D$ in Eq. \eqref{eq:def_D}. 
\begin{theorem}\label{thm:behavior_near_boundary}
Suppose that $\calI \subset \calH $ is an open interval. If $x(h)$ reaches a local maximum at a point $h_0\in \calI$, then $x_0 = x(h_0)$ is a left end point of $\supp(F)$ and $x''(h_0)<0$. Moreover, if $x_0$ is away from $D^\complement$, for sufficiently small $\varepsilon>0$, the density function of $F$ behaves like the square-root function as $f(x) = Q(\sqrt{x-x_0})$ on $(x_0, x_0 + \varepsilon)$ where $Q(x)$ is a real analytic function near $zero$, $Q(0)=0$, and
\[Q'(0) =\frac{1}{\pi}\sqrt{\frac{-2}{x''(h_0)}}\int\frac{d\calB(\tau)}{(\tau - h_0)^2}.   \]

Assume now that $x(h)$ reaches a local minimum at a point $h_0\in \calI$, then $x_0 = x(h_0)$ is a right end point of $\supp(F)$ and $x''(h_0)>0$. Moreover, if $x_0$ is away from $D^\complement$, for sufficiently small $\varepsilon>0$, $f(x) = Q(\sqrt{x_0-x})$ on $(x_0 -\varepsilon , x_0)$ where $Q(x)$ is a real analytic function near $zero$, $Q(0)=0$, and
\[Q'(0) =\frac{1}{\pi}\sqrt{\frac{2}{x''(h_0)}}\int\frac{d\calB(\tau)}{(\tau - h_0)^2}.   \]

\end{theorem}

\begin{proof}
In this proof, we only consider the case where $x(h_0)$ is a local maximum. The case where $x(h_0)$ is a local minimum is treated similarly. We follow the argument in the proof of Theorem 5.2 of \citet{silverstein1995analysis} and Theorem 3.3 of \cite{couillet2014analysis}. If $x(h_0)$ is a local maximum, we prove that $x''(h_0)<0$. When $h$ is near $h_0$, define 
\begin{align*}
P_k(h) &= \gamma \int\frac{u^{k}d\calA(u)}{(1+um_{\calB}(h))^k}, \quad k=1,2,\dots, \\
Q_k(h) &= \int\frac{d\calB(\tau)}{(\tau - h)^k}, \quad k=1,2,\dots.
\end{align*}
It is easy to check that  $P'_k(h) = -k P_{k+1}(h) Q_2(h)$  and $Q'_k(h) = k Q_{k+1}(h)$, $k=1,2,\dots$. With the notation, 
\begin{align*}
&x(h) = h+ P_1(h),\\
&x'(h) = 1 - P_2(h)Q_2(h),\\
&x''(h) = 2 P_3(h) Q_2^2(h) - 2 P_2(h) Q_3(h),\\
&x'''(h) = -6 P_4(h)Q_2^3(h) + 12 P_3(h)Q_2(h)Q_3(h) - 6 P_2(h)Q_4(h).
\end{align*}
Since $x(h_0)$ is a local maximum, $x'(h_0) =  1- P_2(h_0)Q_2(h_0) = 0$. We assume that  $x''(h_0) = 0$ and raise a contradiction. If $x''(h_0) = 0$, $P_3(h_0)Q_2^2(h_0)  = P_2(h_0)Q_3(h_0)$. Then, 
\begin{align*} 
x'''(h_0) &= 6\Big(P_3(h_0) Q_2(h_0)Q_3(h_0) - P_2(h_0) Q_4(h_0)\Big) +  6\Big( P_3(h_0) Q_2(h_0)Q_3(h_0) - P_4(h_0)Q_2^3(h_0) \Big)  \\
& = \frac{6}{Q_2^2(h_0)}\Big(Q_3^2(h_0) - Q_4(h_0)Q_2(h_0)\Big) + \frac{6}{P_2^4(h_0)} \Big(P_3^2(h_0) - P_4(h_0)P_2(h_0)\Big).
\end{align*}
By Cauchy-Schwarz, 
\[ Q_3^2(h_0) - Q_4(h_0)Q_2(h_0) = \left(\int\frac{d\calB(\tau)}{(\tau - h_0)^3} \right)^2 - \int\frac{d\calB(\tau)}{(\tau - h_0)^4} \int\frac{d\calB(\tau)}{(\tau -h_0)^2}\leq 0,\]
with equality holds only when $\calB$ is a Dirac measure. Similarly, 
\[  P_3^2(h_0) - P_4(h_0)P_2(h_0) = \left(\int\frac{u^3d\calA(\tau)}{(1+um_\calB(h_0))^3}\right)^2 - \int\frac{u^4d\calA(\tau)}{(1+um_\calB(h_0))^4}\int\frac{u^2d\calA(\tau)}{(1+um_\calB(h_0))^2}\leq 0,\]
with equality holds only when $\calA$ is a Dirac measure. It follows that $x'''(h_0)<0$ when $\calB$ is not a Dirac measure. Therefore, $x(h_0)$ is an inflection point instead of a local maximum, which is a contradiction. Therefore, we conclude that $x''(h_0)<0$.  

To show that $x_0$ is a left end point of $\supp(F)$, first of all, $x_0$ cannot be an isolated point of $\supp(F)$ since $h_0 \notin \supp(\calB)$. Secondly, by Theorem \ref{thm:increasing_outside_necessary_sufficient}, $(x_0- \varepsilon ,x_0)\subset \supp^\complement(F)$ for some $\varepsilon >0$ since $x'(h) >0$ when $h$ is to the left of $h_0$. If $x_0 \in \supp(F)$, we can immediately conclude that $x_0$ is a left end point of $\supp(F)$. Assume $x_0 \in \supp^\complement(F)$. Again by Theorem \ref{thm:increasing_outside_necessary_sufficient}, we can find $\tilde{h}\in\calH$ such that $x(\tilde{h}) = x_0$ and $x'(\tilde{h})>0$. Using Proposition \ref{prop:x_h_increasing}, it can only happen when $\tilde{h} > h_0$. Indeed, if not, we can find $\varepsilon' >0$ such that $ \tilde{h} < \tilde{h}+ \varepsilon' < h_0 - \varepsilon <h_0$, $x(\tilde{h} + \varepsilon) > x(\tilde{h}) = x(h_0) > x(h_0 - \varepsilon)$, and $x'(\tilde{h} + \varepsilon) >0$, $x'(h_0 - \varepsilon)>0$.  However, Proposition \ref{prop:x_h_increasing} indicates that $x(\tilde{h}  + \varepsilon) < x(h_0  -\varepsilon)$, which is a contradiction.

Since $\calH$ is open and $x'(h)$ is continuous, it follows that we can find $h_*$ such that $h_0 < h_*< \tilde{h}$, $x'(h_*)>0$, and $x(h_*) < x(\tilde{h}) = x_0$. While $x(h)$ is increasing when $h$ is to the left of $h_0$, we can find $h_{**}<h_0$, $x'(h_{**})>0$ but $x(h_{**})> x(h_*)$. However, Proposition \ref{prop:x_h_increasing} indicates that if $h_{**}<h_*$, $x'(h_{**})>0$, $x'(h_*)>0$, we have $x(h_{**}) <x({h_*})$, which is a contradiction. The assumption that $x_0 \in \supp^\complement(F)$ does not hold.

To show the expression of $f(x)$ for $x\in(x_0, x_0+\varepsilon)$ and $x$ is away from $D^\complement$, we can extend $x(h)$ analytically to a neighborhood $\frakN$ of $h_0$ in $\mathbb{C}$. We write this extension as 
\[ z(h) = h+ \gamma \int\frac{ud\calA(u)}{1+um_{\calB}(h)}, \quad h \in \frakN\subset \mathbb{C}.\]
Clearly, $z(h) = x(h)$ when $h\in\frakN$ is real.  On the other hand,  assume $\varepsilon>0$ is sufficiently small so that $(x_0, x_0+\varepsilon) \subset\supp(F)$. For $x\in(x_0, x_0+\varepsilon)$, using Eq. \eqref{eq:main_eq}, we can find $h_{\mathbb{C}}(x)\in\mathbb{C}^+$ with positive imaginary part such that 
\begin{align*}
&\underline{m}(x) = m_{\calB}(h_{\mathbb{C}}(x)) = \int\frac{d\calB(\tau)}{\tau - h_{\mathbb{C}}(x)},\\
&h_{\mathbb{C}}(x) = x - \gamma \int\frac{ud\calA(u)}{1+u\underline{m}(x)}.
\end{align*}
We assume $\varepsilon$ is sufficiently small so that $h_{\mathbb{C}}(x) \in \frakN$. We therefore have $z(h_{\mathbb{C}}(x)) = x$.  Moreover, 
$z'(h_0) = x'(h_0) = 0 $ and $z''(h_0) = x''(h_0)<0$.

Let $\phi(h)$ be an analytic function satisfying $\phi^2(h) =z(h) - x_0$. Then, $\phi(h_0) = 0$ and $(\phi'(h_0))^2 = z''(h_0)/2 =x''(h_0)/2$. While there are two functions satisfying the requirements, we choose $\phi$ such that 
\[\phi'(h_0) = -i\sqrt{ - z''(h_0)/2}.\]
We therefore have 
\[ z(h_{\mathbb{C}}(x))  = x_0 + \phi^2(h_{\mathbb{C}}(x)) = x .\] 
Considering the local analytic inverse $\Phi$ of $\phi$ in a neighborhood of $h_0$, we get that  $h_{\mathbb{C}}(x) =\Phi(\sqrt{x- x_0})$.
Note that the choice of $\phi$ ensures that $\Im h_{\mathbb{C}}(x) = \Im \Phi(\sqrt{x-x_0}) >0$. Moreover, $\Phi$ satisfies $\Phi(0) = h_0$ and $\Phi'(0) = 1/\phi'(h_0) = i (-2/x''(h_0))^{1/2}$.

In a neighborhood of zero, define
\[ Q(z) = \frac{1}{\pi} \Im \int\frac{d\calB(\tau)}{\tau - \Phi(z)}.\]
Clearly, $Q$ is an analytic function near zero as $\Phi(z)$ is near $h_0$ when $z$ is near zero and $h_0\in \calH  \subset \supp^\complement(\calB)$. Then,
\[Q(0) = \frac{1}{\pi} \Im \int\frac{d\calB(\tau)}{\tau - h_0} = 0,\]
\[Q'(0)=\frac{1}{\pi}\sqrt{\frac{-2}{x''(h_0)}}\int\frac{d\calB(\tau)}{(\tau - h_0)^2}. \]
The expression of $f(x)$ follows from the fact that
\[ f(x) = \frac{1}{\pi} \Im \underline{m}(x) = \frac{1}{\pi} \Im\int\frac{d\calB(\tau)}{\tau -  \Phi(\sqrt{x-x_0}) } = Q(\sqrt{x-x_0}).\]
Here, the existence of $f(x)$ when $x\in D$ is ensured by Theorem \ref{thm:main}.

\end{proof}

\appendix
\section{Technical lemmas}\label{appendix:lemmas}

We collect a series of lemmas that are broadly used in the Random Matrix Theory literature. 
\begin{theorem}[{\citet[Theorem 2.1]{silverstein1995analysis}}]\label{thm:Silverstein1995_1}
Let $F$ be a distribution function and $x \in \mathbb{R}$. Suppose $\lim_{z \in \mathbb{C}^{+} \rightarrow x} \Im m_F(z)$ exists. Then, $F$ is differentiable at $x$, and its derivative is $\frac{1}{\pi} \Im m_F(x)$.
\end{theorem}

\begin{theorem}[{\citet[Theorem 2.2]{silverstein1995analysis}}]\label{thm:Silverstein1995_2}
Let $\Omega$ be an open and bounded subset of $\mathbb{R}^n$, let $\Theta$ be an open and bounded subset of $\mathbb{R}^m$, and let $f: \bar{\Omega} \rightarrow \Theta$ be a function, continuous on $\Omega$ ($\bar{\Omega}$ denoting the closure of $\Omega)$. If, for all $x_0 \in \partial \Omega, \lim _{x \in \Omega \rightarrow x_0} f(x)=f\left(x_0\right)$, then $f$ is continuous on all of $\bar{\Omega}$.
\end{theorem}

The following lemma provides a useful bound on the Stieltjes transform when $z$ is away from the real line. 
\begin{lemma}\label{lemma:bound_m_by_Imz}
Suppose that $m(z)$, $z\in\mathbb{C}^+$ is the Stieltjes transform of a probability measure $\mu$, then
\[|m(z)| \leq \frac{1}{\Im z}.\] 
\end{lemma}
\begin{lemma}\label{lemma:point_mass_formula}
Supppose that $\mu$ is a probability measure. The mass of $\mu$ at $x\in\mathbb{R}$ can be retrieved using its Stieltjes transform as
\[\mu(\{x\}) = \lim_{y\downarrow 0} -iy m_{\mu}(x+iy).\]
\end{lemma}
Lemma \ref{lemma:point_mass_formula} can be easily proved using the dominated convergence theorem. See for example, the proof of Proposition 2.2 of \cite{couillet2014analysis} or Theorem 3.1 of \cite{dozier2007analysis}.
\begin{lemma}[Weyl's inequality]\label{lemma:weyl}
Let $A$ and $B$ be two Hermitian matrices. Denote the $j$-th largest eigenvalue of a Hermitian matrix $C$ to be $\tau_j(C)$. Then,
\[ \tau_{i+j-1}(A+B) \leq \tau_i(A)+\tau_j(B) \leq \tau_{i+j-n}(A+B).\] 
\end{lemma}

\section{Additional details in the proof of Lemma \ref{lemma:same_limit_converg}}\label{appendix:additional_proof_lemma_same_limit}

Recall the proof of Lemma \ref{lemma:same_limit_converg} presented in Section \ref{sec:existence_subsec:proof_density}. In this section, we show that $\liminf_{n\to\infty} A_2>0$. 
\begin{align}
& A_2/\gamma=\nonumber\\
&\int \frac{ u^2\left|(u\tilde{m}_n-um_n^*)\tau +(um_n^*z^*_n - u\tilde{m}_n \tilde{z}_n) + (z^*_n - \tilde{z}_n) +(1+u\tilde{m}_n)\gamma\displaystyle\int\frac{wd\calA(w)}{1+w\tilde{m}_n} - (1+u{m}_n^*)\gamma\displaystyle\int\frac{wd\calA(w)}{1+wm_n^*}\right|^2}{\left|(1+um_n^*)(1+u\tilde{m}_n) \Big(\tau -z^*_n + \gamma\displaystyle\int\frac{wd\calA(w)}{1+wm_n^*}\Big)\Big(\tau -\tilde{z}_n + \gamma\displaystyle\int\frac{wd\calA(w)}{1+w\tilde{m}_n}\Big)\right|^2 }\nonumber\\
&d\calB(\tau)d\calA(u)\nonumber\\
& \geq \int \frac{u^2 |(u\tilde{m}_n-um_n^*)\tau +(um_n^*z^*_n - u\tilde{m}_n \tilde{z}_n) + (z^*_n - \tilde{z}_n) +(1+u\tilde{m}_n)\gamma\displaystyle\int\frac{wd\calA(w)}{1+w\tilde{m}_n} - (1+u{m}_n^*)\gamma\displaystyle\int\frac{wd\calA(w)}{1+wm_n^*}|^2}{(1+M|u|)^4 (2\tau^2 + 2K^2)^2}\nonumber\\
&d\calB(\tau)d\calA(u),\nonumber
\end{align}
where $M$ is a finite upper bound on the moduli of $m_n^*$ and $\tilde{m}_n$ and $K$ is a finite upper bound on the moduli of $-z^*_n + \gamma\displaystyle\int\frac{wd\calA(w)}{1+wm_n^*}$ and $-\tilde{z}_n + \gamma\displaystyle\int\frac{wd\calA(w)}{1+w\tilde{m}_n}$ (See Lemma \ref{lemma:boundedness_beta}).
By Fatou's lemma, 
\begin{align}
\liminf_{n\to\infty} A_2/\gamma
&\geq\int \frac{ \left|(\underline{\tilde{m}}-\underline{m}^*)(u\tau -ux) +(1+u\underline{\tilde{m}})\gamma\displaystyle\int\frac{wd\calA(w)}{1+w\underline{\tilde{m}}} - (1+u\underline{m}^*)\gamma\displaystyle\int\frac{wd\calA(w)}{1+w\underline{m}^*}\right|^2}{u^{-2} (1+M|u|)^4 (2\tau^2 + 2K^2)^2}d\calB(\tau)d\calA(u),\label{eq:long_eq}
\end{align}
 For the numerator on the right-hand side, we can obtain a matrix form as 
\begin{align*}
(\underline{\tilde{m}}-\underline{m}^*)(u\tau -ux) +(1+u\underline{\tilde{m}})\gamma\int\frac{wd\calA(w)}{1+w\tilde{m}_n} - (1+u\underline{m}^*)\gamma\int\frac{wd\calA(w)}{1+w\underline{m}^*} = [\begin{matrix} \tau -x, &1\end{matrix}] \Delta \left[\begin{matrix} 1\\ u\end{matrix}\right],
\end{align*}
where 
\begin{align*}
\Delta = \left[\begin{matrix}0 && & \underline{\tilde{m}}-\underline{m}^* \\[10pt] \gamma\displaystyle\int\frac{wd\calA(w)}{1+w\underline{\tilde{m}}} -\gamma\displaystyle\int\frac{wd\calA(w)}{1+w\underline{m}^*}& & &\underline{\tilde{m}}\gamma\displaystyle\int\frac{wd\calA(w)}{1+w\underline{\tilde{m}}} - \underline{m}^*\gamma\displaystyle\int\frac{wd\calA(w)}{1+w\underline{{m}^*}}  \end{matrix}\right].
\end{align*}
All together
\begin{align*}
\mbox{Right-hand side of Eq. \eqref{eq:long_eq}} = \int \operatorname{Tr}\{H(\tau) \Delta G(u) \Delta^*\} d\calB(\tau)d\calA(u) = \Big[\operatorname{Tr} \int H(\tau)d\calB(\tau)\Big] \Delta \Big[\int G(u)d\calA(u) \Big] \Delta^*,
\end{align*}
where 
\[H(\tau) = \frac{1}{4(\tau^2 + K^2)^2}\left[\begin{matrix} (\tau-x)^2 & \tau -x \\ \tau-x & 1 \end{matrix} \right], \]
\[G(u) = \frac{|u|^2}{(1+M|u|)^4} \left[\begin{matrix} 1 & &&u \\ u&& & u^2\end{matrix} \right].\]
Recall that a $2\times 2$ matrix is positive definite if and only if at least one diagonal element is positive and the determinant is positive.
Using H\"{o}lder's inequality
\[\det\left(\int H(\tau) d\calB(\tau)\right) = \int \frac{(\tau-x)^2}{4(\tau^2+K^2)^2}d\calB(\tau) \int \frac{1}{4(\tau^2+K^2)^2}d\calB(\tau) - \left(\int \frac{\tau-x}{4(\tau^2+K^2)^2}d\calB(\tau)\right)^2 \]
is strictly positive, given that $\calB(\tau)$ is not a Dirac measure.  It follows that $\int H(\tau) d\calB(\tau)$ is positive definite. By the same argument, we also have that $\int G(u)d\calA(u)$ is positive definite. While $\Delta$ is non-zero due to the assumption that $\underline{\tilde{m}} - \underline{m}^*\neq 0$, we have the right hand side of Eq. \eqref{eq:long_eq} is strictly positive. The proof is complete.

\section{Proof of Theorem \ref{thm:existence_of_masses}}\label{appendix:proof_theorem_existence_masses}

In this section, we complete the proof of Theorem \ref{thm:existence_of_masses} and show the sufficiency of the condition $\calB(\{b\}) - \gamma(1-\calA(\{0\}))$. It is more convenient to take the perspective of matrix algebra. We shall use Weyl's inequality to bound the eigenvalues of $F^{\bW_n}$ near $b$. Suppose we have $\calB(\{b\}) - \gamma(1-\calA(\{0\})) >0$. Denote the e.d.f of the eigenvalues of $\bW_n$, $\bA_n$, $\bB_n$ and $n^{-1}\bX_n^*\bA_n\bX_n$ to be $F^{\bW_n}$, $F^{\bA_n}$, $F^{\bB_n}$ and $F^{n^{-1}\bX^*_n\bA_n\bX_n}$ respetively. 
In view that $F$, $\calA$, $\calB$ is the weak limit of $F^{\bW_n}$, $F^{\bA_n}$, $F^{\bB_n}$ almost surely (see Introduction).  Moreover, since $n^{-1}\bX_n^*\bA_n\bX_n$ can be view as a special case of $\bB_n + n^{-1}\bX_n^*\bA_n\bX_n$ when $\bB_n=0$, we have that $F^{n^{-1}\bX^*_n\bA_n\bX_n}$ converges almost surely to a probability function denoted as $\calG$ weakly. 

We first assume $\bA_n$ is nonnegative definite and therefore $\calA$ is supported on $\mathbb{R}^+\cup\{0\}$. Then, by Proposition 2.2 of \citet{couillet2014analysis}, $\calG(\{0\}) = 1- \min\{1, \gamma(1-\calA(\{0\}))\} = 1- \gamma(1-\calA(\{0\}))$, since $\gamma(1-\calA(\{0\}))<\calB(\{b\})\leq 1$. 
Consider a sufficiently small constant $\epsilon>0$ such that $\calG$ is continuous at $\epsilon/2$, $F$ and $\calB$ are continuous at $b\pm \epsilon$ and $b\pm \epsilon/2$. Note that $F$ and $\calA$ can have at most countably many discontinuity points. 
Let $K = K(n,\epsilon)$ be the number of eigenvalues of $n^{-1}\bX_n^*\bA_n\bX_n$ in $[0,\epsilon/2]$ and let $M = M(n,\epsilon)$ be the number of eigenvalues of $\bB_n$ on $(b-\epsilon, b+\epsilon/2]$. Then,
\begin{align*}
&\frac{K}{n} =  F^{n^{-1}\bX_n^*\bA_n\bX_n}(\epsilon/2) \longrightarrow \calG(\epsilon/2), \quad \mbox{almost surely},\\
&\frac{M}{n} =  F^{\bB_n}(b+\epsilon/2) - F^{\bB_n}(b-\epsilon) \longrightarrow \calB(b+\epsilon/2) - \calB(b-\epsilon), \quad \mbox{almost surely}.
\end{align*}
Given that $\calB(\{b\}) - \gamma (1-\calA(\{0\})) = \calB(\{b\}) - 1 +\calG(\{0\})>0$, let $\kappa = (\calB(\{b\}) - 1 +\calG(\{0\}))/2$. For all sufficiently small $\epsilon$ we have $\calB(b+\epsilon/2) - \calB(b-\epsilon) - 1 +\calG(\epsilon/2) >\kappa$. Therefore,
\[ P\left( \frac{M -n +K}{n} >\kappa\right) \longrightarrow 1.\] 
Now, consider the eigenvalues of $\bW_n = \bB_n + n^{-1}\bX^*_n \bA_n \bX_n$ near $b$. Using Weyl's inequality (Lemma \ref{lemma:weyl}), we get that at least $M-(n-K)$ eigenvalues of $\bW_n$ are in $(b-\epsilon, b+\epsilon]$. It follows that  for all sufficiently small $\epsilon$,
\[P\left(F^{\bW_n}(b+\epsilon)  - F^{\bW_n}(b-\epsilon) > \kappa\right) \longrightarrow1,\]
and therefore $F(b+\epsilon) - F(b-\epsilon) >\kappa$. We therefore have $F(\{b\})\geq \kappa$. 

If $\bA_n$ has negative eigenvalues, say $u_1^n \geq \cdots u_{p_1}^{n} \geq 0 > u_{p_1+1}^n \geq \cdots \geq u_p^{n}$. Define $\bA_{(1)} = \operatorname{diag}(u_1^n,\dots u_{p_1}^n)$ and $\bA_{(2)} = \operatorname{diag}(u_{p_1+1}^n, \cdots, u_{p}^n)$. Then
\[ n^{-1} \bX_n^* \bA_n\bX_n = n^{-1} \bX_{(1)}^* \bA_{(1)}\bX_{(1)} + n^{-1}\bX_{(2)}^* \bA_{(2)}\bX_{(2)},\]
where $\bX_{(1)}$ is the first $p_1$ rows of $\bX_n$ and $\bX_{(2)}$ is the last $p-p_1$ rows of $\bX_n$. Note that $\bX_{(1)}^* \bA_{(1)}\bX_{(1)}$ and $\bX_{(2)}^* \bA_{(2)}\bX_{(2)}$ are independent. Let $\tilde{\bB}_n  = \bB_n + n^{-1}\bX_{(1)}^* \bA_{(1)}\bX_{(1)}$. Then, the previous results apply to $\tilde{\bB}_n$ and $-\bW_n= -\tilde{\bB}_n + n^{-1}\bX_{(2)}^*(-\bA_{(2)})\bX_{(2)}$. We omit the details.

\bibliographystyle{apalike}
\bibliography{draft}

\end{document}